\documentclass{article}
\usepackage{arxiv}

\usepackage[utf8]{inputenc}

\usepackage{lineno}
\usepackage{standalone}
\usepackage{hyperref}
\hypersetup{
	colorlinks,
	linkcolor={black},
	citecolor={black},
	urlcolor={blue!80!black}
}

\usepackage{caption}
\usepackage{subcaption}
\usepackage{stmaryrd}
\usepackage{amssymb,amsmath,amstext}
\usepackage{graphicx}
\usepackage{booktabs}

\usepackage{tikz}
\usetikzlibrary{decorations.pathmorphing,decorations.pathreplacing,decorations.markings,patterns,calc}
\usetikzlibrary{shapes.geometric, shapes.arrows, arrows}
\usetikzlibrary{spy,backgrounds}
\usetikzlibrary{external}
\usetikzlibrary{intersections}
\usetikzlibrary{tikzmark}

\tikzstyle{startstop} = [rectangle, rounded corners, minimum width=3cm, minimum height=1cm,text centered, draw=black, fill=red!30]
\tikzstyle{io} = [trapezium, trapezium left angle=70, trapezium right angle=110, minimum width=3cm, minimum height=1cm, text centered, draw=black, fill=blue!30]
\tikzstyle{process} = [rectangle, minimum width=3cm, minimum height=1cm, text centered, draw=black, fill=orange!30]
\tikzstyle{decision} = [diamond, minimum width=3cm, minimum height=1cm, text centered, draw=black, fill=green!30]
\tikzstyle{arrow} = [->,>=stealth]

\usepackage{import}
\usepackage{pdfpages}
\usepackage{transparent}
\usepackage{xcolor}

\usepackage{textcomp}
\usepackage{siunitx}
\usepackage[most]{tcolorbox}

\usepackage{bm}
\usepackage[]{datatool}
\usepackage{calc}
\usepackage{textgreek}
\usepackage{bigints}
\usepackage{physics}
\usepackage{interval}
\usepackage{accents}
\usepackage{mathtools,eqparbox}
\usepackage{indentfirst}
\usepackage[algo2e,ruled,linesnumbered,vlined,commentsnumbered,longend]{algorithm2e}
\usepackage{enumerate}

\usepackage{pgfplots}
\pgfplotsset{compat=newest}
\usepgfplotslibrary{groupplots,dateplot}
\usepgfplotslibrary{colormaps}
\usepgfplotslibrary{colorbrewer}
\usepgfplotslibrary{polar}
\usepgfplotslibrary{fillbetween}

\usepackage[absolute,overlay]{textpos}

\usepackage{tikz-dimline}

\usepackage[autostyle]{csquotes}

\usepackage[colorinlistoftodos]{todonotes}

\DeclareMathOperator*{\argmin}{arg\,min}

\DeclareMathOperator{\sgn}{sgn}
\DeclareMathOperator{\sign}{sign}

\DeclareMathOperator{\diag}{diag}

\newcommand{\tn}{\mathrm{n}} 
\newcommand{\ts}{\mathrm{s}} 

\newcommand{\bv}{\boldsymbol{v}}

\newcommand{\bu}{\boldsymbol{u}} 
\newcommand{\bn}{\boldsymbol{n}} 
\newcommand{\bx}{\boldsymbol{x}} 
\newcommand{\ubx}{\mathbf{x}} 
\newcommand{\by}{\mathbf{y}} 
\newcommand{\bY}{\mathbf{Y}} 
\newcommand{\uby}{\mathbf{y}} 
\newcommand{\bz}{\mathbf{z}} 

\newcommand{\bSigma}{\boldsymbol{\Sigma}} 

\newcommand{\ba}{\boldsymbol{a}}

\newcommand{\bepsilon}{\boldsymbol{\varepsilon}} 

\newcommand{\bb}{\boldsymbol{b}}

\newcommand{\jump}[1]{\left\llbracket #1 \right\rrbracket} 
\newcommand{\bphi}{\boldsymbol{\varphi}}

\newcommand{\rom}[1]{\uppercase\expandafter{\romannumeral #1\relax}}

\intervalconfig {
	soft open fences ,
}

\usepackage{amsthm}

\theoremstyle{plain}

\theoremstyle{definition}
\newtheorem{definition}{Definition}[section]

\theoremstyle{remark}
\newtheorem{remark}{Remark}

\pgfplotsset{
	legend image with text/.style={
			legend image code/.code={%
					\node[anchor=center] at (0.3cm,0cm) {#1};
				}
		},
}

\usepackage{doi}

\usepackage[
	natbib=true,
	backend=biber,
	giveninits=true,
	date=year,
	doi=false,
	url=false,
	isbn=false,
	eprint=false,
	bibstyle=ieee,
	citestyle=numeric-comp,
]{biblatex}
\title{Robust design optimization for enhancing delamination resistance of composites}

\author{
  Sukhminder Singh$^{1, 2}$\thanks{Corresponding author, email address: \texttt{sukhminder.singh@fau.de}}\ , Lukas Pflug$^1$, Julia Mergheim$^3$, and Michael Stingl$^2$ \\ \\
$^1$Competence Unit for Scientific Computing (CSC), \\
Friedrich-Alexander-Universit\"at~Erlangen-N\"urnberg~(FAU) \\ \\
$^2$Chair of Applied Mathematics (Continuous Optimization), \\
Friedrich-Alexander-Universit\"at~Erlangen-N\"urnberg~(FAU) \\ \\
$^3$Institute of Applied Mechanics, \\
Friedrich-Alexander-Universit\"at~Erlangen-N\"urnberg~(FAU)
}

\begin{document}









\maketitle
\begin{abstract}
	Recent developments in the field of computational modeling of fracture have opened up possibilities for designing structures against failure.
	A special case, called interfacial fracture or delamination, can occur in loaded composite structures where two or more materials are bonded together at comparatively weak interfaces.
	Due to the potential crack growth along these interfaces, the structural problem suffers from snap-back/snap-through instabilities and bifurcations with respect to the model parameters, leading to noisy and discontinuous responses.
	For such a case, the design optimization problem for a selected quantity of interest is ill-posed, since small variations in the design parameters can lead to large jumps in the structural response.
	To this end, this paper presents a stochastic optimization approach to maximize delamination resistance that is less sensitive to small perturbations of the design and thereby leads to a robust solution.
	To overcome the intractability of Monte Carlo methods for estimating the expected value of the expensive-to-evaluate response function, a global, piecewise-constant surrogate is constructed based on nearest-neighbor interpolation that is iteratively refined during the optimization run.
	We found that by taking a large stochastic region at the beginning of the optimization and gradually reducing it to the desired one can help overcome poor local optima.
	Our results demonstrate the effectiveness of the proposed framework using an example of shape optimization of hard inclusions embedded in a double-cantilever beam, which significantly enhances delamination resistance.
\end{abstract}


\section{Introduction}\label{sec:introduction}

Structural defects in the form of sharp cracks tend to grow in highly stressed regions and in those regions that provide easy pathways for crack propagation, for instance, adhesive interfaces in laminated composites.
In the latter case, introducing architected heterogeneities or voids at the micro- or meso-scale along the interfaces is expected to retard or even stop the growth of the crack, thereby increasing the resistance to interfacial fracture or \emph{delamination}.
The size, shape and material properties of the heterogeneities can be chosen based on the solution of a structural optimization problem in which the objective is to maximize or minimize a \emph{quantity of interest (QoI)} that characterizes the fracture resistance of the interfaces.

On the downside for the formulation of such an optimization problem, structural analysis involving quasi-static crack growth exhibits two major characteristics:
snap-back/snap-through instabilities, and
bifurcations~\cite{wriggers2008nonlinear,NonlinearDynamStroga2018} with respect to the model parameters resulting in a switching of crack patterns (see \cite{StochasticPhasGerasi2020,MaterialOptimiSingh2021} for examples).
Such behavior is also evident in other structural problems, e.g., related to nonlinear buckling~\cite{StabilityBifuPignat1998}, crash~\cite{SurrogateModelBoursi2018}, or damage~\cite{StabilityOfStBazant2010}.
As a result, design optimization problem formulations without some special treatment of the above structural phenomena, are \emph{discontinuous}.
Here, the discontinuities refer to the large jumps in the objective function value under small perturbations of the design parameters, rendering the optimization problem ill-posed~\cite{IllPosedProblBaier1994}.

To overcome the challenge of ill-posedness, attempts have been made to simplify the structural analysis model, and use density- (e.g., SIMP~\cite{GeneratingOptiBendso1988}) and level-set-based~\cite{ALevelSetMetAllair2002} topology optimization techniques to obtain enhanced fracture-resistant structures.
Historically, this has been accomplished with great success by incorporating stress constraints~\cite{StudyOnTopoloCheng1992,TopOptStressConstraintsBendsoe1998,StressBasedToLeCh2010,AnEnhancedAggLuoY2013,StressConstraiAmir2017,StressBasedShPicell2018}, improving the resistance to crack nucleation in areas of high stress concentration.
Similar techniques have been developed for stationary cracks to improve their resistance to grow further~\cite{OptimizationOfGuGr2016,TopologyOptimiKang2017,MinimizingCracKlarbr2018,FractureStrengHuJi2020,OnTailoringFrZhang2022}.
Due to the restriction of crack evolution, such optimization problems do not face the challenges of structural instabilities and bifurcations.
Moreover, the analytical sensitivities of the objective with respect to the design parameters can be reliably exploited by efficient gradient-based optimization methods.

Over the past decade, numerous studies have attempted to incorporate propagating cracks into the framework of topology optimization, primarily by regularizing the description of discrete cracks using a scalar phase field~\cite{RevisitingBritFrancf1998,AReviewOnPhaAmbati2015}.
Typically, the objective function involves a structural response function integrated over a finite number of fixed load steps, e.g.,
integrated mechanical work~\cite{TopologyOptimiDaDa2018,TopologyOptimiDaDa2020,LevelSetTopolWuCh2020,APathDependenWuCh2021},
integrated fracture surface energy~\cite{ANovelTopologRuss2020}, or
integrated elastic energy~\cite{TopologyOptimiDesai2022}.
For the computation of reliable analytical sensitivities of the objective with respect to the design, these studies (except in Desai~et~al~\cite{TopologyOptimiDesai2022}) used a staggered algorithm~\cite{APhaseFieldMMiehe2010} together with a history variable field for elastic energy to robustly solve the fracture problem.
In Desai~et~al~\cite{TopologyOptimiDesai2022}, however, the authors depended on a backtracking algorithm~\cite{TheVariationalBourdi2008} to achieve continuity of the potential energy of the system w.r.t.\ pseudo-time, while trusting the analytical sensitivities derived without considering structural instabilities.
Similarly, in Singh~et~al~\cite{MaterialOptimiSingh2021}, that addresses material optimization for fracture resistance of cohesive interfaces, the authors achieved stable crack growth by locally linearizing the traction-separation law for each loading step, enabling calculation of exact gradients of the objective function with respect to the design.
As shown in the same work, a major drawback of this approach is that it requires very small loading steps to resolve fast evolving cracks, and furthermore the method leads to spurious local optima in the objective function landscape close to the design points for which the structural system exhibits bifurcation.

\begin{figure}
	\centering
	\def\svgwidth{1\columnwidth}
	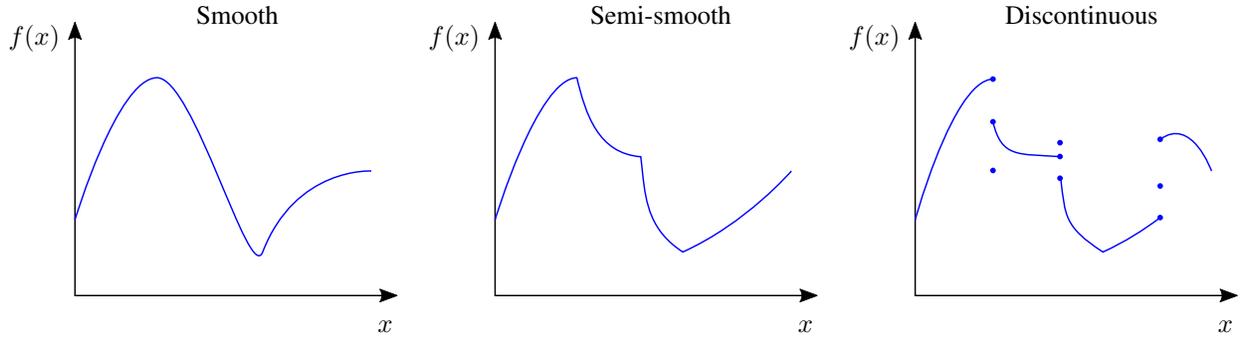

	\caption{
		Illustration of different classes of structural response $f$ in terms of a design parameter $x$ --- smooth, semi-smooth and discontinuous, in the absence of numerical noise.
		The structural responses for the delamination problem exhibiting bifurcation phenomenon belong to the third category.
		At the locations of discontinuities, multiple solutions to the delamination problem may exist, leading to distinct responses.
	}
	\label{fig:function_types}
\end{figure}

\begin{figure}
	\begin{center}
		\includegraphics[width=0.4\textwidth]{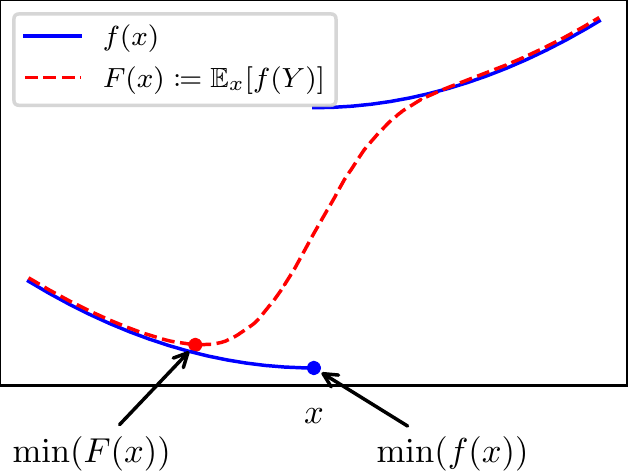}
	\end{center}
	\caption{
	A schematic showing the minimum points of a discontinuous response surface $f(x)$ and its smoothed counterpart $F(x) \coloneqq \mathbb{E}_x[f(Y)]$.
	The minimum of $f(x)$ lies at the point of discontinuity which leads to high sensitivity of the response w.r.t.\ $x$, whereas the minimum of $F(x)$ is located away from the discontinuity, which leads to robustness.
	}
	\label{fig:robust_solution_schematic}
\end{figure}

To design a structure with optimized resistance to crack propagation, it would be better not to incorporate a regularization which alters the physics of the system and also the solution should be robust in the presence of discontinuities.
Fig.~\ref{fig:function_types} shows three different classes of cost functions encountered in structural optimization problems: \emph{smooth}, \emph{semi-smooth} and \emph{discontinuous}.
For the smooth and semi-smooth variants, well established gradient- and surrogate-based optimization methods exist in the literature that help finding a local or the global minimum of the cost function.
However, for the discontinuous case, it is challenging for both types of methods to solve the optimization problem with a reasonable computational expense.
In addition, the solution found may not be robust in the sense that it may lie at or near the discontinuities, making the cost value very sensitive to small variations in the parameters of the structural problem.
Given environmental uncertainties, it is imperative to find solutions which are robust against deviations in the design parameters distributed according to some probability density function.
For these reasons,
we shift from deterministic to stochastic optimization approaches.
In the framework of stochastic optimization, instead of finding a minimum of a scalar cost function $f(\mathbf{x})$, we aim to minimize its expectation over some random variable~$\mathbf{Y}$:
\begin{equation}
	\min_\mathbf{x}\ \left\{F(\mathbf{x})
	\coloneqq \mathbb{E}_\mathbf{x}[f(\mathbf{Y})] \right\},
	\label{eq:expectation}
\end{equation}
which may be viewed as a \emph{stochastic relaxation} of the original deterministic optimization problem (see Fig.~\ref{fig:robust_solution_schematic}).

A challenge of the stochastic optimization approach is that it requires many evaluations of the QoI, for instance, for the evaluation of the high-dimensional integral over the probability distribution of the stochastic variables.
Monte Carlo techniques are the standard approach due to the fact that these methods are straightforward to implement and do not fall under the \emph{curse of dimensionality}.
Unfortunately, the number of samples required to estimate the mean of the QoI and its gradient with respect to the design variables can be prohibitive for direct use of a gradient based optimization algorithm.
Improvements such as using the adaptive multiple importance sampling technique to incorporate all the samples generated during an optimization run~\cite{ADerivativeFrMaggia2018} can greatly reduce the computational costs, but the number of model evaluations often remains prohibitively high.
Efforts have been made to reduce this computational burden by replacing the expensive-to-evaluate response function by global approximations, called \emph{surrogate models} or \emph{metamodels}.
A surrogate model substitutes the true functional relationship $f(\mathbf{x})$ with a mathematical expression $\tilde{f}(\mathbf{x})$ generated by a few sample points ${\{\mathbf{x}_i, f(\mathbf{x}_i)\}}_{i=1}^N$ and is much cheaper to evaluate.
During the optimization run, the surrogate model is refined with each iteration by adding more samples in the interesting regions in the design space for which the computer experiments can be run independently and in parallel.
To name a few, common surrogate models used for optimization include
\emph{polynomial regression model},
\emph{Gaussian process (GP) regression} (or \emph{Kriging})~\cite{ComparisonOfRSimpso1998,FacilitatingPrKoch2002,KrigingAssisteRaponi2019},
\emph{radial basis function (RBF)}~\cite{RationalRadialJakobs2009},
and so on.
However, the effectiveness of the surrogate models diminishes significantly when the true response surface is discontinuous with respect to the design parameters, and the discontinuities can lead to oscillations due to Gibbs phenomenon.
Discontinuity detection strategies such as Voronoi-based~\cite{VpsVoronoiPiRushdi2017} or support vector machines (SVM)~\cite{RobustUncertaiWildey2019} based domain-decomposition try to find the piecewise-continuous subdomains, enabling fitting of separate continuous surrogate models depending upon the number of discontinuities.

In this paper, we simply use a piecewise-constant surrogate based on nearest-neighbor interpolation~\cite{TheUniformConDevroy1978}, constructed from normally distributed points sampled randomly around the design points encountered in an optimization run.
The main advantage lies in the elimination of Gibbs oscillations in the presence of discontinuities in the original function.
Using specialized data structures, e.g.\ $k$-d tree~\cite{AnalysisApproxNNMount1999}, nearest-neighbor interpolation also enables fast updating and evaluation of the model.
In a nutshell, we try to combine the good attributes from these two methods: the ability of stochastic relaxation to provide robust optimized designs and the ability of nearest-neighbor interpolation to approximate discontinuous functions.
To the best of our knowledge, none of the previous works has taken into account the robustness of the final design for optimized fracture resistance.
The idea is to apply a gradient-based optimization algorithm on the smoothed, stochastically relaxed version of the piecewise-constant approximation of the response surface.
To overcome spurious local optima due to numerical noise, an adaptive stochastic relaxation technique, similar to Covariance Matrix Adaptation Evolution Strategy (CMA-ES)~\cite{TheCmaEvolutiHansen2016}, is used with an adaptively evolving trust region as in the Successive Response Surface Method (SRSM)~\cite{OnTheRobustneStande2002}.
Although adaptive sampling techniques can be employed to significantly improve the convergence rate of surrogates, especially in low-dimensional problems, their adoption and comparison with random sampling in the context of optimization are beyond the scope of this work.

The remainder of the paper is organized as follows.
We begin with Section~\ref{sec:delamination_model} explaining the mathematical formulation of the delamination problem along with an appropriate metric for delamination resistance.
Section~\ref{sec:robust_design_optimization} describes the formulation of the stochastic optimization problem and an algorithm to find a robust local optimum.
Section~\ref{sec:numerical_experiments} presents some numerical experiments, testing the optimization method on a heterogeneous double-cantilever beam problem with two, six, and twelve design degrees of freedom (DoFs).
Finally, in Section~\ref{sec:conclusion}, we conclude the paper highlighting the advantages and limitations of the optimization model.

\section{Mathematical model for simulating interfacial fracture}%
\label{sec:delamination_model}
In this section, we briefly describe the main concepts for modeling interfacial fracture in a heterogeneous material domain.
We use an energy minimization-based approach to evaluate the evolution of cracks along the material interfaces, wherein the interfacial tractions are determined by an exponential cohesive law.
We also briefly discuss the spatial discretization scheme based on the eXtended Finite Element Method (XFEM)~\cite{TheExtendedGeFries2010} to model displacement discontinuities along the interfaces, and a technique to solve the discretized system.
Subsequently, we define a convenient QoI for structural optimization, which is chosen as a metric for delamination resistance.

\subsection{The energy approach}%
\begin{figure}
	\centering
	\def\svgwidth{0.8\columnwidth}
	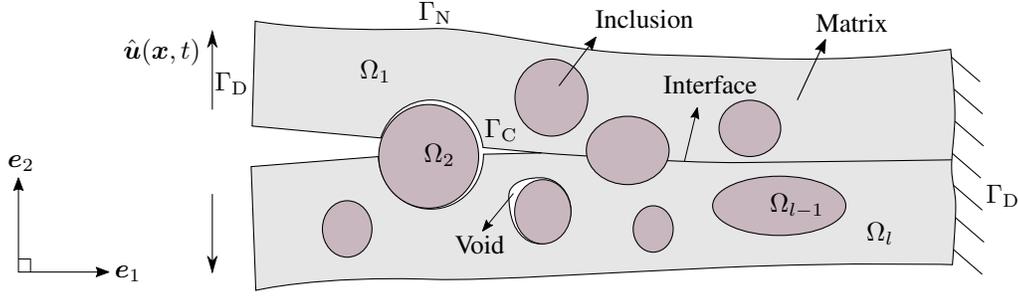

	\caption{Schematic of a continuum domain $\Omega$ with cohesive interfaces $\Gamma_\mathrm{C}$.}
	\label{fig:delamination_schematic}
\end{figure}

\label{sec:energy_approach}
Let $\Omega \subset \mathbb{R}^2$ be an open and bounded domain representing the configuration of a composite body in two dimensions, as shown in Fig.~\ref{fig:delamination_schematic}.
We assume that $\Omega$ is divided into $L \in \mathbb{N}$ homogeneous, linear-elastic subdomains $\Omega_l, l = 1,\dots,L$, interwoven by adhesive interfaces denoted by the set $\Gamma_\mathrm{C} \coloneqq \cup_{l=1}^L \partial \Omega_l \setminus \partial \Omega$.
The domain boundary $\partial \Omega$ is divided into two disjoint sets, $\Gamma_\mathrm{D}$ and $\Gamma_\mathrm{N}$, representing Dirichlet and Neumann boundaries, respectively.
The state of this mechanical system for a given set of boundary conditions is given by the displacement field $\bu : \overline{\Omega} \rightarrow \mathbb{R}^2$, which admits discontinuities on $\Gamma_\mathrm{C}$.
The interfacial opening at a point $\bx \in \Gamma_\mathrm{C}$ with unit vector $\bn(\bx)$ pointing in the direction perpendicular to the interface is given by the displacement jump $\jump{\bu(\bx)} \coloneqq \bu(\bx^+) - \bu(\bx^-)$, where
$\bx^\pm = \lim_{\epsilon \rightarrow 0^+} \bx \pm \epsilon \bn(\bx)$.
On the Dirichlet boundary, a time-dependent displacement field $\hat{\bu} : \Gamma_\mathrm{D} \times [0, T] \rightarrow \mathbb{R}^2, \ T \in \mathbb{R}_{> 0}$ is prescribed, whereas for the sake of simplicity of the model, we assume that no external body and surface forces act on the system.

To define the state of the system at a given (pseudo-) time instant $t \in [0, T]$, we search for solutions to the mechanical problem in the space
\begin{align}
	\mathcal{S} \coloneqq \left\{ \bv \in L^2(\Omega; \mathbb{R}^2) : \bv \vert_{\Omega_l} \in H^1(\Omega_l; \mathbb{R}^2),\ l=1,\dots,L \right\}
	\label{eq:state_solution_space}
\end{align}
by minimizing the total potential energy functional $\mathcal{E} : \mathcal{S} \rightarrow \mathbb{R}_{\geq 0}$, given by
\begin{align}
	\mathcal{E}(\bu) \coloneqq \sum_{l = 1}^L \int_{\Omega_l} \Psi_l(\bepsilon(\bu(\bx))) \, \dd{\bx}
	+ \int_{\Gamma_\mathrm{C}} \mathcal{G}(\jump{\bu(\bx)}) \, \dd{s},
	\label{eq:total_potential_energy}
\end{align}
with $\bu \equiv \hat{\bu}(\,\cdot\,, t)$ on $\Gamma_\mathrm{D}$.
Here, $\Psi_l$ represents the stored elastic energy density, $\mathcal{G}$ is the cohesive potential function, and $\dd{s}$ denotes the surface measure over $\Gamma_\mathrm{C}$.
The elastic energy density at a point $\bx \in \Omega_l,\ l=1,\dots,L$ is given by
\begin{equation}
	\Psi_l(\bepsilon)
	\coloneqq \frac{1}{2} \bepsilon : \mathbb{C}_l : \bepsilon,
	\label{eq:stored_elastic_energy}
\end{equation}
where $\bepsilon \coloneqq \nabla^\mathrm{sym} \bu$ is the second-order, symmetric strain tensor and $\mathbb{C}_l$ is the fourth-order elasticity tensor.
Assuming isotropic behavior of the bulk, the elasticity tensor is defined in terms of Lam\'e parameters $\lambda_l$ and $\mu_l$ as
\begin{equation}
	\mathbb{C}_l \coloneqq \lambda_l \boldsymbol{I} \otimes \boldsymbol{I} + 2\mu_l \boldsymbol{1}^\mathrm{sym}.
	\label{eq:linear_elasticity_tensor}
\end{equation}
Here, $\boldsymbol{I}$ is the second-order identity tensor and $\boldsymbol{1}^\mathrm{sym}$ is a fourth-order unit tensor, defined as $1^\mathrm{sym}_{ijkl} = [\delta_{ik}\delta_{jl} + \delta_{il}\delta_{jk}]/2$, where $\delta_{ij}$ is the Kronecker delta with the property $\delta_{ij} = 1$ if $i=j$, and $\delta_{ij} = 0$ if $i \neq j$.

To define the constitutive behavior of the interfaces, we follow Xu and Needleman~\cite{VoidNucleationXuX1993} and adopt an exponential cohesive potential as a function of interfacial separation, expressed as
\begin{equation}
	\begin{aligned}
		\mathcal{G}(\jump{\bu}) & \coloneqq \phi(\delta_\tn \coloneqq \jump{\bu} \cdot \bn, \delta_\ts \coloneqq \jump{\bu} \cdot \bn^\perp)                                                                                                                                                               \\
		                        & \coloneqq \phi_\tn + \phi_\tn \exp(-\frac{\delta_\tn}{\delta_\tn^*})  \left\{ \left[1-r+\frac{\delta_\tn}{\delta_\tn^*}\right] \frac{1-q}{r-1} - \left[q+ \frac{r-q}{r-1} \frac{\delta_\tn}{\delta^*_\tn}\right] \exp(-\frac{\delta_\ts^2}{{\delta_\ts^*}^2}) \right\} ,
		\label{eq:needleman_fracture_energy}
	\end{aligned}
\end{equation}
where $\delta_\mathrm{n}$ and $\delta_\mathrm{s}$ are the displacement jumps in the directions normal and tangential to the interface, respectively.
The unit vector $\bn^\perp$ is perpendicular to $\bn$ and points in the direction parallel to the interface.
The parameters $r$ and $q \coloneqq \phi_\ts/\phi_\tn$ are the coupling parameters, where $\phi_\mathrm{n}$ is the work of normal separation and $\phi_\mathrm{s}$ is the work of tangential separation.
The parameters $\delta_\tn^*$ and $\delta_\ts^*$ denote the critical interfacial openings.
For details about the behavior of this constitutive model and the resulting cohesive traction given by $\nabla \phi$, please refer to \cite{VoidNucleationXuX1993}.

To simulate delamination of the structure under a given displacement controlled loading, we first discretize time over $N_T$ time steps, given by $t^k = k \Delta t, \ k=1,\dots,N_T$, where $\Delta t \coloneqq T / N_T$ is the time-step size.
Denoting the state of the system by displacement field $\bu^k$ at time step $k$,
we write the structural problem as:
\begin{definition}[Incremental Potential Energy Minimization]
	\label{def:energy_minimization}
	Given a prescribed boundary displacement field $\hat{\bu}(\cdot, t^k)$ on $\Gamma_\mathrm{D}$, find
	\begin{equation}
		\begin{aligned}
			\bu^k                \in & \argmin_{\bv \in \mathcal{S}}\ \mathcal{E}(\bv),                                           \\
			                         & \text{such that } \bv(\bx) = \hat{\bu}(\bx , t^k) \quad \forall \bx \in \Gamma_\mathrm{D}.
		\end{aligned}
	\end{equation}
\end{definition}


\subsection{Spatial discretization and solution method}
\label{sec:spatial_discretization}

As indicated earlier, we model the displacement discontinuities at the interfaces using XFEM, which is a flexible numerical approach for general interface problems.
For completeness, we introduce here only its basic formulation.

\begin{figure}
	\centering
	\def\svgwidth{1\columnwidth}
	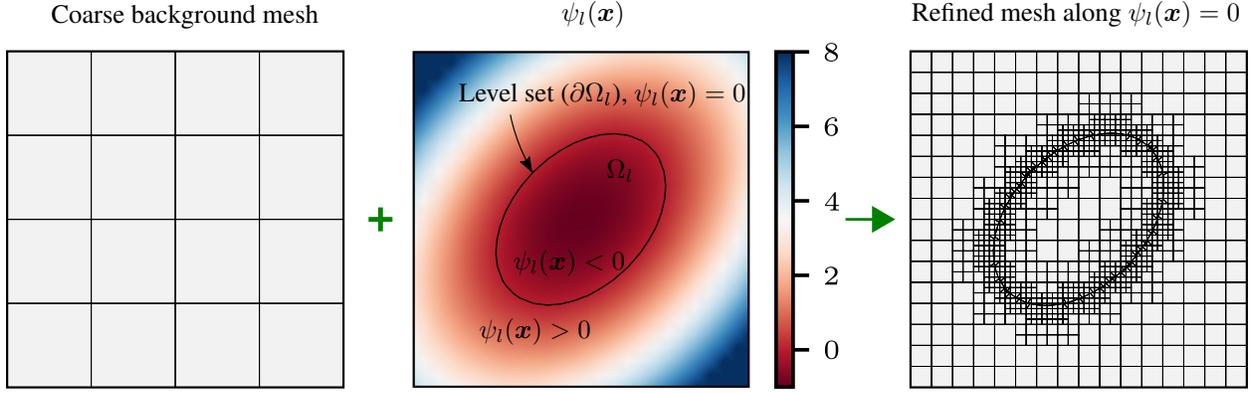

	\caption{
		Illustration of spatial discretization of a square domain with an elliptical heterogeneity $\Omega_l$ using a level-set function.
		(Left) structured coarse background mesh,
		(center) level-set function defining the interface between the matrix and the heterogeneity,
		and (right) adaptively refined mesh along the interface.
		(For more information about the color references in this figure, the reader is referred to the digital version of this article.)
	}
	\label{fig:spatial_discretization_schematic}
\end{figure}

\begin{figure}
	\centering
	\def\svgwidth{0.8\columnwidth}
	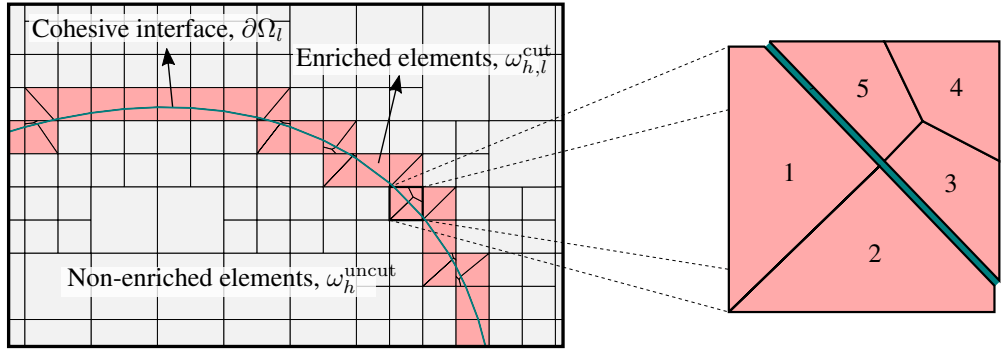

	\caption{(Left) Enriched elements allowing jump discontinuity along the bonded interface,
		and (right) subdivision of a cut cell into \num{5} quadrilaterals to apply a standard quadrature rule for integration over the sub-elements~\cite{OnImplementationXFEMCarraro2015}.}
	\label{fig:enriched_elements}
\end{figure}

To define the subdomains $\Omega_l,\ l=1,\dots,L$, we use smooth level-set functions $\psi_l : \Omega \rightarrow \mathbb{R}$, where $\psi_l = 0$ on the respective interfaces.
Starting from a coarse background mesh (see Fig.~\ref{fig:spatial_discretization_schematic}), the finite-elements intersected by the interfaces, referred to as \emph{cut cells}, are successively refined to a predefined depth $H$ or to a level where each of the finite elements is cut by only one interface.
Let $\omega_{h,l}^\mathrm{cut} \subset \Omega$ be a subdomain covered by all the cells cut by an interface represented by the level-set function $\psi_l$, and $\omega_{h}^\mathrm{uncut} \coloneqq \Omega \setminus \cup_{l=1}^{L} \omega_{h,l}^\mathrm{cut}$ be the part of the spatial domain covered by all the non-intersected cells, as shown in Fig.~\ref{fig:enriched_elements}.
At time step $k$, we express the discretized displacement field $\bu_h$ at any point $\bx \in \Omega$ as
\begin{equation}
	\begin{aligned}
		\bu^k(\bx) \approx \bu_h(\bx; \ba^k, \bb^k) \coloneqq
		\begin{cases}
			\sum_{i \in I_a}  a_i^k \bphi_i(\bx)                                                      & \text{for } \bx \in \omega_h^\mathrm{uncut},   \\
			\sum_{i \in I_a}  a_i^k \bphi_i(\bx) + \sum_{i \in I_b} b_i^k \xi_{l,i}(\bx) \bphi_i(\bx) & \text{for } \bx \in \omega_{h,l}^\mathrm{cut},
		\end{cases}
	\end{aligned}
	\label{eq:discretized_displacement_field}
\end{equation}
where $\bphi_i$ are the conventional finite element shape functions.
The field vectors $\ba^k \in \mathbb{R}^{|I_a|}$ and $\bb^k \in \mathbb{R}^{|I_b|}$ represent conventional and enriched nodal displacements, respectively, where $I_a$ and $I_b$ are the corresponding index sets.
To model the discontinuities within the cut cells, a usual XFEM approach uses an enrichment function of the form
\begin{equation}
	\xi_{l,i}(\bx) \coloneqq \sgn(\psi_l (\bx)) - \sgn(\psi_l (\boldsymbol{X}_i)),
	\label{eq:discontinuous_enrichment_function}
\end{equation}
where $\boldsymbol{X}_i$ is the position of the finite element node corresponding to the $i^\mathrm{th}$ degree of freedom.
For more details on the implementation of XFEM for interface problems, please refer to~\cite{OnImplementationXFEMCarraro2015}.

Denoting by $I_\mathrm{D}$ the index set of constrained degrees of freedom corresponding to the prescribed displacements at the Dirichlet boundary, the energy minimization problem in Definition~\ref{def:energy_minimization} is restated in spatially discretized form as below:
\begin{definition}[Incremental Potential Energy Minimization for Spatially Discretized System]%
	\label{def:energy_minimization_discretized}
	Given nodal displacements
	$\hat{a}^k_i,\ i \in I_\mathrm{D}$ at time-step $k$,
	find $\bu_h^k \equiv \bu_h(\,\cdot\,; \ba^k, \bb^k)$ such that
	\begin{equation}
		\begin{aligned}
			\{\ba^k, \bb^k\} \in & \argmin_{\ba, \bb}\ \mathcal{E}(\bu_h(\,\cdot\,; \ba, \bb)),           \\
			                     & \text{such that}\ a_i = \hat{a}_i^k\quad \forall i \in I_{\mathrm{D}}.
		\end{aligned}
		\label{eq:discretized_potential_energy_minimization}
	\end{equation}
\end{definition}

Since the energy functional $\mathcal{E}$ is twice-differentiable with respect to the unknowns $\ba$ and $\bb$, the above optimization problem can be solved using a continuous, gradient-based optimization solver.
For the numerical studies presented in this paper, we used \texttt{IPOPT}~\cite{OnTheImplemenWachte2006}, which is an open-source software for solving large scale nonlinear optimization problems.
To achieve faster convergence, the solution $\{\ba^{k-1}, \bb^{k-1}\}$ from the previous time step is used as an initial guess, providing a warm start to the optimization procedure.
The gradient $\nabla_{\{\boldsymbol{a}, \boldsymbol{b}\}} \mathcal{E}(\bu_h(\,\cdot\,; \ba, \bb))$ and the Hessian $\nabla^2_{\{\ba, \bb\}}\mathcal{E}(\bu_h(\,\cdot\,; \ba, \bb))$ of the energy functional, as required by the optimizer, can be derived manually, or one can employ an automatic differentiation (AD) package, for example, \texttt{Sacado}~\cite{AutomaticDiffeBartle2006}, within a finite-element assembler code.

\begin{remark}
	The energy minimization problem defined by Definition~\ref{def:energy_minimization_discretized} admits non-unique solutions due to the non-convex nature of the potential energy functional.
	The solution obtained is possibly one of the many local minima, determined by the internal parameters of the state problem solver (here, \texttt{IPOPT}).
	In order to have a deterministic solution, i.e.\ same solution for multiple problem solves, these parameters must be kept constant throughout a design optimization run.
\end{remark}

\subsection{Quantity of interest}%
\label{sec:quantity_of_interest}
The aim of this study is to optimize heterogeneous structures for enhanced resistance against interfacial fracture.
For a structure subjected to a displacement-controlled loading, the resistance to its interfacial fracture can be characterized by the amount of external mechanical work that must be expended to displace a part of the boundary by a given amount.
The mechanical work of the applied displacement is given by its scalar product with the reaction force exerted by the respective boundary.
For the purposes of structural optimization presented in this paper, we consider the QoI to be the negative of the external mechanical work, numerically integrating external mechanical power over time using the trapezoidal rule~\cite{davis2007methods}, as
\begin{align}
	\mathcal{W}(\bu_h^0,\dots,\bu_h^{N_T}) \coloneqq - \frac{1}{2}  \sum_{k=1}^{N_T} \sum_{l=1}^L \int_{\Gamma_\mathrm{D} \cap \partial \Omega_l} \left[\boldsymbol{\tau}_l(\bu_h^k(\bx)) + \boldsymbol{\tau}_l(\bu_h^{k-1}(\bx))\right] \cdot \left[\bu_h^k(\bx) - \bu_h^{k-1}(\bx)\right] \, \dd{s},
	\label{eq:mechanical_work}
\end{align}
where $\boldsymbol{\tau}_l(\bu) \coloneqq [\mathbb{C}_l : \bepsilon(\bu)] \cdot \bn$ is the surface traction on $\Gamma_\mathrm{D} \cap \partial \Omega_l$ with unit normal vector $\bn$.
The negative sign implies that we want to maximize the mechanical work by minimizing $\mathcal{W}$ and thus improving the delamination resistance.
In the following sections, the terms QoI and structural response are used interchangeably.

\section{Robust design optimization}
\label{sec:robust_design_optimization}
The delamination problem introduced in Section~\ref{sec:delamination_model} assumes that the parameters involved in the model are deterministic, i.e., known for certain, or produced exactly to value.
These include parameters related to the geometry and material of the structure, as well as the discretization and internal parameters of the solution algorithm.
Obviously, this approach is an idealization of the real-life scenario, where uncertainties are inherent.
To analyze a structural problem with non-unique solution, it is crucial to identify states that the system can assume with high probability given a particular set of model parameters.
In the context of structural optimization, this aspect becomes even more important as considering only one state solution and neglecting the others can lead to a final design that shows a completely different behavior than predicted in a real life situation.
Even if the underlying structural problem has a unique solution for the optimized design, small variations in the parameters can have large effects on the structural response.
On the contrary, a \emph{robust design} is characterized by minimal influence of fluctuations on the system behavior (for example, see~\cite{TopologyOptimiDeSu2020}).

In this section, bearing in mind the characteristics of the delamination problem, we introduce a stochastic optimization problem formulation.
We explain the nearest-neighbor interpolation scheme to construct a piecewise-constant surrogate of the structural response.
To solve the stochastic optimization problem, a surrogate-based algorithm is presented along with the search gradients to be used for the minimization of the objective.

\subsection{Formulation of the stochastic optimization problem}
Let us assume that $\ubx \in {\mathcal{B}}$ be a vector of $d \in \mathbb{N}$ design variables (not to be confused with the spatial variable $\bx$ used in the previous sections), where
\begin{equation}
	\mathcal{B} \coloneqq \left\{ \ubx \in \mathbb{R}^d \mid x_i^\mathrm{min} \leq x_i \leq x_i^\mathrm{max}, i=1,\dots,d \right\}
	\label{eq:feasible_design_set}
\end{equation}
represents simple box constraints, with design bounds $\mathbf{x}^\mathrm{min}, \mathbf{x}^\mathrm{max} \in \mathbb{R}^d$.
To formulate a stochastic optimization problem, we introduce $\mathbf{Y} \in \mathbb{R}^d$ as a vector of random (stochastic) variables representing a perturbed design around $\mathbf{x}$, following a distribution with some additional parameters $\boldsymbol{\Theta}$.
The probability density function for the distribution is given by $\rho(\,\cdot\,; \mathbf{x}, \boldsymbol{\Theta}) : \mathbb{R}^d \rightarrow \mathbb{R}_{\geq 0}$.
If $f(\uby)$ denotes a deterministic model mapping parameters to a QoI, for example $\mathcal{W}$ in Eq.~\eqref{eq:mechanical_work}, and assuming $\rho(\mathbf{y};\,\cdot\,,\,\cdot\,) = 0$ for an infeasible design $\mathbf{y} \in \mathbb{R}^d$, we write the constrained stochastic optimization problem as
\begin{equation}
	\min_{\mathbf{x} \in \mathcal{B}}\
	\Bigg\{ F(\ubx)
	\coloneqq \mathbb{E}_{\mathbf{x}, \boldsymbol{\Theta}}[f(\mathbf{Y})]
	=	\int_{{\mathbb{R}^d}} f(\by) \rho(\by; \ubx, \boldsymbol{\Theta})\, \dd{\by}\Bigg\}.
	\label{eq:min_F}
\end{equation}

A number of options for choosing the probability distribution (e.g., Weibull, Poisson, normal, lognormal, exponential, etc.) exist which can be suitable for particular optimization problems.
One of the most important and simplest ones is the \emph{multivariate normal} (or \emph{Gaussian}) distribution defined by the probability density
\begin{align}
	p(\by; \ubx, \bSigma) \coloneqq \frac{1}{{(2\pi)}^{\frac{d}{2}} {\vert \bSigma \vert}^{\frac{1}{2}}} \exp(-\frac{1}{2} {[\by-\ubx]}^\top \bSigma^{-1} [\by-\ubx]),
	\label{eq:multivariate_normal_pdf}
\end{align}
with mean $\ubx$ and symmetric, positive-definite covariance matrix $\bSigma \in \mathbb{R}^{d \times d}$.
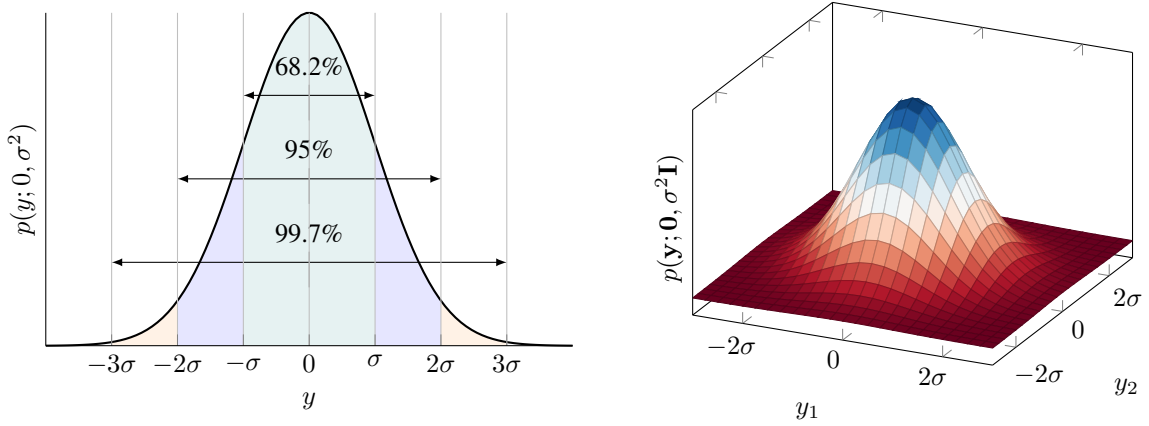
\begin{figure}
	\centering
	\begin{subfigure}{0.52\textwidth}
		\begin{center}
			\pgfmathdeclarefunction{gauss}{2}{\pgfmathparse{1/(#2*sqrt(2*pi))*exp(-((x-#1)^2)/(2*#2^2))}%
}
\begin{tikzpicture}[draw=black, text=black]
	\begin{axis}[no markers, domain=0:10, samples=100,
			axis lines*=left, xlabel={$y$}, ylabel={$p(y; 0, \sigma^2)$},
			width=\textwidth,
			height=0.7\textwidth,
			xtick={-3, -2, -1, 0, 1, 2, 3}, ytick=\empty,
			xticklabels={$-3\sigma$, $-2\sigma$, $-\sigma$, $0$, $\sigma$, $2\sigma$, $3\sigma$},
			enlargelimits=false, clip=false, axis on top,
			grid = major]
		\addplot [fill=teal!10, draw=none, domain=-3:3] {gauss(0,1)} \closedcycle;
		\addplot [fill=orange!10, draw=none, domain=-3:-2] {gauss(0,1)} \closedcycle;
		\addplot [fill=orange!10, draw=none, domain=2:3] {gauss(0,1)} \closedcycle;
		\addplot [fill=blue!10, draw=none, domain=-2:-1] {gauss(0,1)} \closedcycle;
		\addplot [fill=blue!10, draw=none, domain=1:2] {gauss(0,1)} \closedcycle;
		\addplot [thick, domain=-4:4] {gauss(0,1)};
		\addplot[latex-latex] coordinates {(-1,0.3) (1,0.3)};
		\addplot[latex-latex] coordinates {(-2,0.2) (2,0.2)};
		\addplot[latex-latex] coordinates {(-3,0.1) (3,0.1)};
		\node[coordinate, pin={68.2\%}] at (axis cs: 0, 0.27){};
		\node[coordinate, pin={95\%}] at (axis cs: 0, 0.17){};
		\node[coordinate, pin={99.7\%}] at (axis cs: 0, 0.07){};
	\end{axis}
\end{tikzpicture}
		\end{center}
	\end{subfigure}
	\begin{subfigure}{0.45\textwidth}
		\begin{center}
			\def\centerx{0}
\def\centery{0}
\begin{tikzpicture}[draw=black, text=black]
	\begin{axis}[
			width=\textwidth,
			xlabel=$y_1$,
			ylabel=$y_2$,
			zlabel={$p(\mathbf{y}; \mathbf{0}, \sigma^2 \mathbf{I})$},
			xtick={-4, -2, 0, 2, 4},
			ytick={-4, -2, 0, 2, 4},
			ztick=\empty,
			zticklabels ={},
			xticklabels={$-4\sigma$, $-2\sigma$, $0$, $2\sigma$, $4\sigma$},
			yticklabels={$-4\sigma$, $-2\sigma$, $0$, $2\sigma$, $4\sigma$},
			colormap/RdBu-11,
		]
		\addplot3[surf,domain=-3:3,domain y=-3:3]
		{1/(2*pi)* exp(-( (x-\centerx)^2 + (y-\centery)^2)/2 )};
	\end{axis}
\end{tikzpicture}
		\end{center}
	\end{subfigure}
	\caption{
		(Left) Univariate normal probability density function $p(y; x, \sigma^2) = \frac{1}{\sqrt{2\pi}\sigma} \exp(-\frac{{(y-x)}^2}{2\sigma^2})$ with mean $x=0$ and standard deviation $\sigma$ where \SI{99.7}{\percent} of the probability mass lies within the $\pm 3 \sigma$ range.
		(Right) Multivariate normal probability density function~\eqref{eq:multivariate_normal_pdf} in two dimensions with mean $\mathbf{x}=\mathbf{0}$ and covariance matrix $\bSigma = \sigma^2 \mathbf{I}$.
	}
	\label{fig:gaussian}
\end{figure}
For the current study, we adopt this distribution and for simplicity, we assume that the random variables $Y_i, i=1,\dots,d$  are independent and identically distributed with a fixed standard deviation $\sigma \in \mathbb{R}_{>0}$, or covariance $\sigma^2 \mathbf{I}$ with identity matrix $\mathbf{I}$.
Here, the standard deviation $\sigma$ measures the magnitude of the random perturbations of each of the design variables.
A well known property of the univariate normal distribution is that, $\SI{99.7}{\percent}$ of the probability mass lies within the $\pm 3\sigma$ range (see Fig.~\ref{fig:gaussian}), which allows us to safely truncate it to this limit without much loss of accuracy.
Bearing in mind the possibility of having infeasible design points in $\mathbb{R}^d$, we define a restricted integration domain for calculating the expectation in Eq.~\eqref{eq:min_F}:
\begin{align}
	\mathcal{B}_\delta & \coloneqq \left\{ \uby \in \mathbb{R}^d \mid x_i^\mathrm{min} - \delta \sigma \leq y_i \leq x_i^\mathrm{max} + \delta \sigma,\ i=1,\dots,d \right\},
\end{align}
where the parameter $\delta \in \mathbb{R}_{> 0}$ is chosen large enough (e.g., $ \delta \geq 3$) such that
\begin{equation}
	\int_{\mathcal{B}_\delta} p(\by; \ubx, \sigma^2 \mathbf{I})\, \dd{\by} \approx 1 \quad \forall \mathbf {x} \in \mathcal{B},
	\label{eq:normalization_approx}
\end{equation}
and small enough ($\delta \leq \infty$) that a perturbed design $\mathbf{y} \in \mathcal{B}_\delta$ around any $\mathbf{x} \in \mathcal{B}$ is still feasible.
This results in a reformulation of the optimization problem, expressed as:
\begin{definition}[Stochastic Optimization Problem]%
	\label{def:optimization_problem}
	Given standard deviation parameter $\sigma \in \mathbb{R}_{>0}$ for design perturbations, find $\ubx \in \mathcal{B}$ such that
	\begin{equation}
		\min_{\mathbf{x} \in \mathcal{B}}\
		\Bigg\{ F(\ubx)
		\coloneqq \mathbb{E}_{\mathbf{x}, \sigma^2 \mathbf{I}}[f(\mathbf{Y})]
		\approx	\int_{{\mathcal{B}_\delta}} f(\by) p(\by; \ubx, \sigma^2 \mathbf{I})\, \dd{\by}\Bigg\} .
		\label{eq:objective_function_truncated_normal}
	\end{equation}
\end{definition}

\begin{figure}
	\centering
		\includegraphics[width=\textwidth]{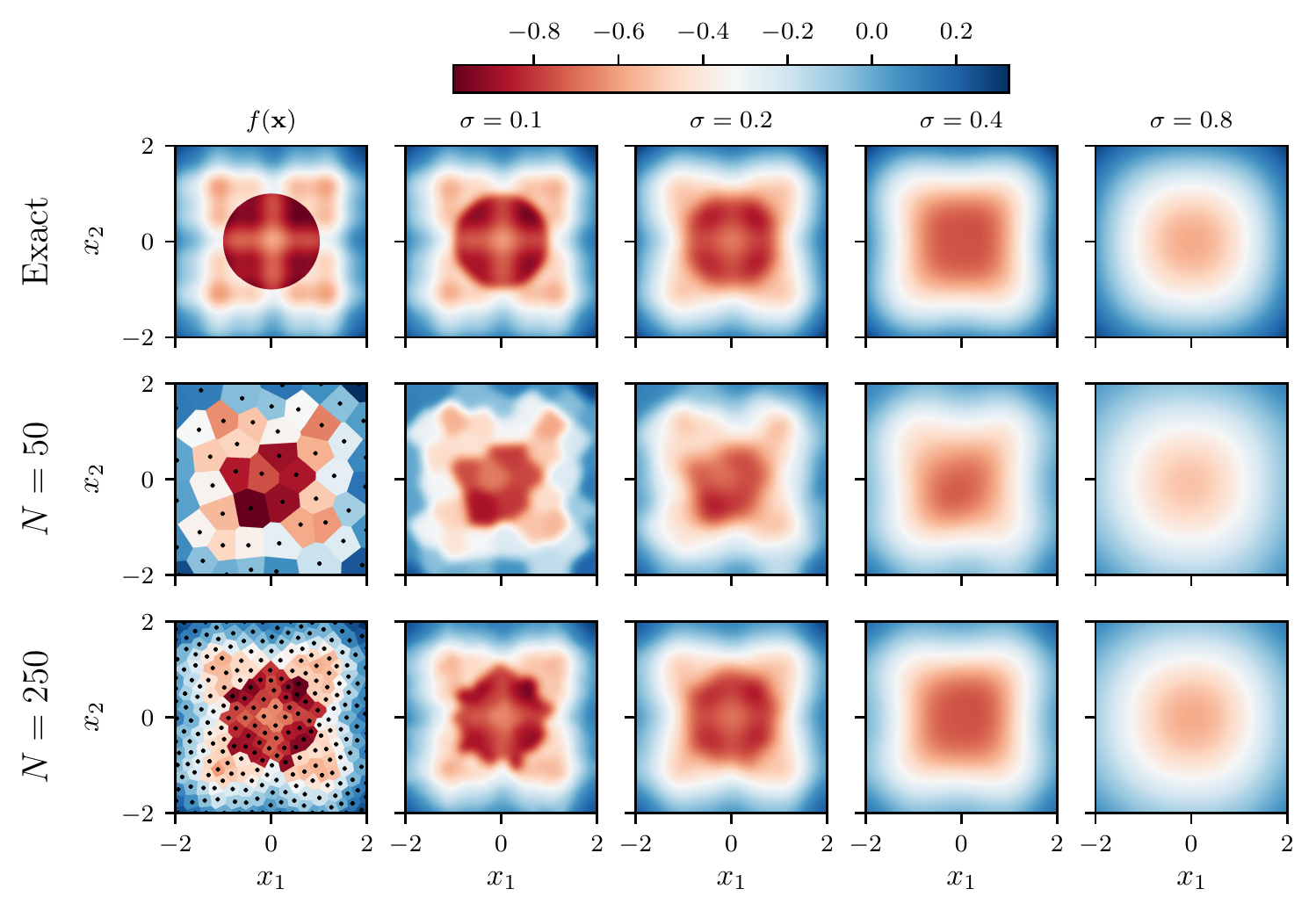}
	\caption{
		Visualizing Herbie-Step function and nearest-neighbor interpolants generated with $N=50$ and $N=250$ well-spaced samples, along with their smoothed counterparts with different values of the smoothing parameter $\sigma$.
		(For more information about the color references in this figure, the reader is referred to the digital version of this article.)
	}
	\label{fig:surrogate_smoothing}
\end{figure}

\begin{remark}[Scaling of design space]
	In the above problem formulation, we have considered a single standard deviation parameter $\sigma$ for all the random variables, which may not be appropriate in every practical situation.
	To account for this, we can define a separate $\tilde{\sigma}_i$ for each of the random variables $Y_i,\ i=1,\dots,d$, leading to a diagonal covariance matrix $\bSigma = \diag(\{\tilde{\sigma}_i^2\}_{i=1}^d)$.
	The stochastic region $\mathcal{B}_\delta$ must then be adjusted accordingly.
	As an alternative, one can scale the design space a-priori such that the stochastic variables $Y_i$ follow a common standard deviation in the optimization framework.
	For the sake of convenience, we take the latter approach in this paper.
\end{remark}

\begin{remark}[Relation to Gaussian smoothing]
	The presented stochastic optimization problem (Definition~\ref{def:optimization_problem}) can be interpreted as \emph{optimization via Gaussian smoothing}~\cite{RandomGradientNester2017,pmlr-v97-maheswaranathan19a,ADerivativeFrMaggia2018}, which is one of the many existing derivative-free optimization (DFO)~\cite{IntroductionToConn2009,DerivativeFreeLarson2019} methods.
	Gaussian smoothing based algorithms first smooth the landscape of the objective function with a $d$-dimensional Gaussian convolution
	and then estimate the gradient of the smoothed objective function through random perturbations of the design.
	As shown in Fig.~\ref{fig:surrogate_smoothing}, using the example of the Herbie function superimposed with a step function, the Gaussian convolution with a wide smoothing radius enables nonlocal exploration, reducing the local minima effect and improving the characterization of global patterns in the objective function landscape.
	In the limit of $\sigma \rightarrow 0$, the smoothed function approaches the true QoI.
	From this perspective, the chosen standard deviation $\sigma$ can also be thought of as a smoothing parameter for the QoI.
	Instead of sequentially scaling $\sigma$ down to $0$, we use it here as a fixed parameter in the context of stochastic optimization.
\end{remark}

\SetKwComment{Comment}{/* }{ */}

\subsection{Optimization algorithm}\label{sec:optimization_algorithm}

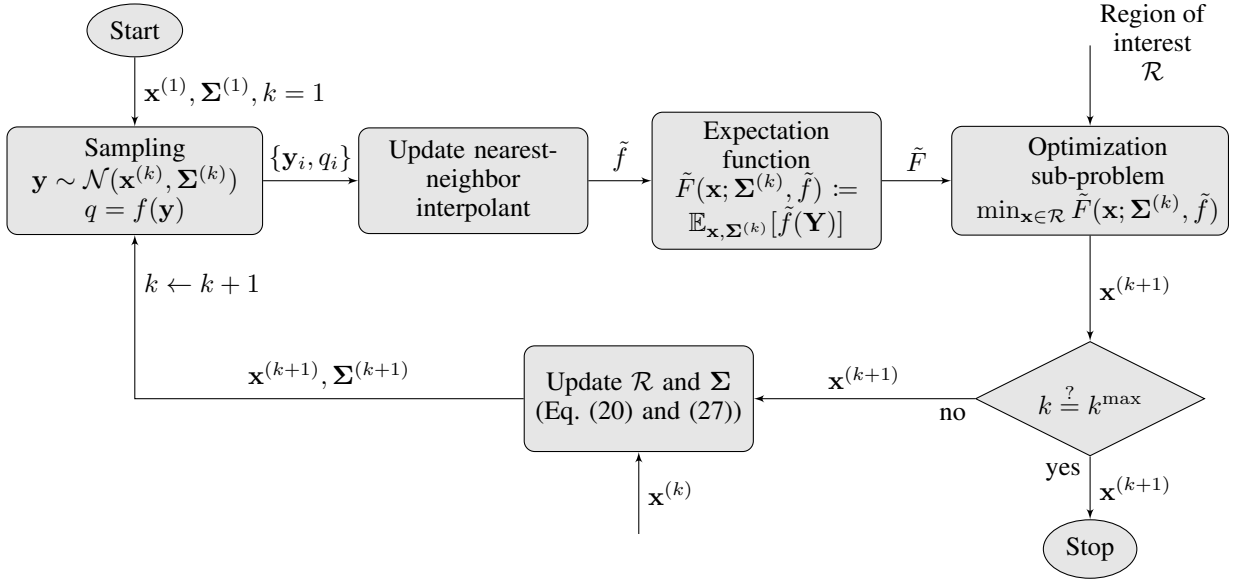
\begin{figure}
	\begin{center}

\tikzstyle{decision} = [diamond, draw, fill=black!10,
node distance=2.9cm, inner sep=4pt, aspect=2]
\tikzstyle{block} = [rectangle, draw, fill=black!10,
text width=8em, text centered, rounded corners, minimum height=4em,
node distance=3.9cm,
]
\tikzstyle{line} = [draw, -latex']
\tikzstyle{cloud} = [draw, ellipse, fill=black!10, node distance=2.9cm,
minimum height=2em]

\begin{tikzpicture}[text=black, draw=black, scale=1, every node/.style={scale=1}]
	\node [cloud] (start) {Start};
	\node [block, below of=start, node distance=2.0cm, text width=9em] (sampling) {\begin{tabular}{c} Sampling                                             \\
			$\mathbf{y} \sim \mathcal{N}(\ubx^{(k)}, \bSigma^{(k)})$ \\ $q = f(\uby)$\end{tabular}};
	\node [block, right of=sampling, node distance=4.5cm] (surrogate) {Update nearest-neighbor interpolant};
	\node [block, right of=surrogate, text width=8em, node distance=3.9cm] (smoothing) {Expectation function $\tilde{F}(\mathbf{x}; \boldsymbol{\Sigma}^{(k)}, \tilde{f}) \coloneqq \mathbb{E}_{\mathbf{x}, \boldsymbol{\Sigma}^{(k)}}[\tilde{f}(\mathbf{Y})]$};
	\coordinate [above of=smoothing, above=0.8cm] (above_smoothing);
	\node [block, right of=smoothing, node distance=4.3cm, text width=9.8em] (optimize) {\begin{tabular}{c} Optimization \\ sub-problem \\ $\min_{\mathbf{x} \in \mathcal{R}} \tilde{F}(\mathbf{x};\bSigma^{(k)},\tilde{f})$ \end{tabular}};
	\coordinate [above of=optimize, above=0.8cm] (above_optimize);

	\node [decision, below of=optimize] (decide) {$k \stackrel{?}{=} k^{\mathrm{max}}$};
	\node [cloud, below of=decide, node distance=2cm] (stop) {Stop};
	\node [block, left of=decide, node distance=6cm] (update) {Update $\mathcal{R}$ and $\bSigma$ (Eq.~\eqref{eq:region_of_interest} and \eqref{eq:update_covariance_matrix})};

	\path [line] (start) -- node[right, midway] {$\mathbf{x}^{(1)}, \bSigma^{(1)}, k=1$} (sampling);
	\path [line] (sampling) -- node [above, midway] {$\{\mathbf{y}_i, q_i\}$} (surrogate);
	\path [line] (surrogate) -- node[above, midway] {$\tilde{f}$} (smoothing);
	\path [line] (smoothing) -- node[above, midway] {$\tilde{F}$} (optimize);

	\path [line] (above_optimize) -- node[right=-0.2cm, pos=0] {\begin{tabular}{c} Region of \\ interest \\ $\mathcal{R}$ \end{tabular}} (optimize);
	\coordinate [below of=update, below=0.8cm] (below_update);
	\path [line] (below_update) -- node[midway, right] {$\ubx^{(k)}$} (update);
	\path [line] (optimize) -- node[midway, right] {$\mathbf{x}^{(k+1)}$} (decide);
	\path [line] (decide) -- node[below, pos=0.1] {no} node[above, midway] {$\ubx^{(k+1)}$} (update);
	\path [line] (update) -| node[pos=0.25, above] {$\ubx^{(k+1)}, \bSigma^{(k+1)}$} node[pos=0.85, right] {$k \gets k+1$} (sampling);
	\path [line] (decide) -- node[left, near start] {yes} node[right, midway] {$\ubx^{(k+1)}$} (stop);

\end{tikzpicture}
	\end{center}
	\caption{
		Flowchart of the optimization procedure. The inputs are the starting design $\ubx^{(1)}$ corresponding to the best Latin-hypercube sample and the covariance matrix $\bSigma^{(1)}$.
	}
	\label{fig:algorithm_flowchart}
\end{figure}

In the following, we explain a surrogate-based algorithm for solving the optimization problem defined by~\eqref{eq:objective_function_truncated_normal}.
The optimization procedure proposed in this work is tailored to efficiently search for a robust design while preventing premature convergence to a poor local optimum in the case of a multi-modal objective landscape.
Fig.~\ref{fig:algorithm_flowchart} gives an overview of the optimization algorithm, whereas a pseudocode is delineated in Algorithm~\ref{alg:optimization_algorithm}.

To find a minimum of the objective $F$ for a given standard deviation $\sigma$, we solve a sequence of stochastically relaxed optimization sub-problems.
To this end, we define a smooth function $\tilde{F} : \mathbb{R}^d \rightarrow \mathbb{R}$,
\begin{align}
	\tilde{F}(\ubx; \bSigma, \tilde{f}) \coloneqq \mathbb{E}_{\mathbf{x}, \bSigma}[\tilde{f}(\mathbf{Y})] = \int_{\mathbb{R}^d} \tilde{f}(\by) p(\by; \ubx, \bSigma)\, \dd{\by} ,
	\label{eq:tilde_F}
\end{align}
where $\tilde{f} : \mathbb{R}^d \rightarrow \mathbb{R}$ is a \emph{nearest-neighbor (NN) interpolant}~\cite{ScatteredDataWendla2004},
built on some data
$\mathcal{D} \coloneqq \{\mathcal{Y}, \mathcal{Q}\}$
with data values $\mathcal{Q} \coloneqq {\{q_i \coloneqq f(\mathbf{y}_i)\}}_{i=1}^N$, evaluated at some scattered data sites $\mathcal{Y} \coloneqq \{\mathbf{y}_1,\dots,\mathbf{y}_N\}$.
For the numerical experiments presented in this paper, distances are measured w.r.t.\ the Euclidean norm to search for nearest neighbors.
Comparison with other distance metrics, e.g.\ Manhattan distance, Mahalanobis distance, etc., is beyond the scope of this work.
Given a sufficiently large number of data points, the nearest-neighbor interpolant $\tilde{f}$ approximates the true response function $f$.

Here, the main difference between $F$ and $\tilde{F}$ is that we have replaced the true response $f$ by its approximation $\tilde{f}$, and instead of fixing the covariance matrix to $\sigma^2 \mathbf{I}$, we have allowed it to remain as a variable.
Having a variable $\bSigma$ allows us to start the optimization with a bigger stochastic region and then successively reduce it to the desired one.
In the case, $\bSigma = \sigma^2 \mathbf{I}$ and $\tilde{f} \approx f$ for a sufficiently large $N$, we recover the true objective function:
\begin{equation}
	F(\ubx) \approx \tilde{F}(\ubx; \bSigma = \sigma^2 \mathbf{I}, \tilde{f} \approx f).
	\label{eq:objective_func_equals_approx}
\end{equation}

\begin{algorithm2e}
	\SetAlgoLined
	\DontPrintSemicolon
	\caption{Stochastic optimization algorithm}\label{alg:optimization_algorithm}
	\KwIn{
		$\mathbf{x}^\mathrm{min}$,
		$\mathbf{x}^\mathrm{max}$,
		$\sigma$,
		$\sigma^\mathrm{max}$,
		$n_0$, 
		$\alpha$, 
		$\delta$, 
		$k^\mathrm{max}$,
		$\gamma^\mathrm{pan}$, 
		$\gamma^\mathrm{osc}$, 
		$\eta$, 
		$\beta$ 
	}
	\KwResult{Optimized design $\mathbf{x}^{*}$}

	$\mathcal{Y} \coloneqq \{\by_i \sim \mathcal{U}(\ubx^\mathrm{min}, \ubx^\mathrm{max}) \}_{i=1}^{n_0}$ \tcp*{Generate uniform random samples using LHS}
	$\mathcal{Q} \coloneqq \{q_i = f(\mathbf{y}_i)\}_{i=1}^{n_0} $ \tcp*{Run delamination simulations and compute QoI}
	$\mathbf{x}^{(1)} \gets \mathbf{y} \in \mathcal{Y}\ \text{ s.t. } f(\mathbf{y}) = \min(\mathcal{Q}) $ \tcp*{Initialize starting design}
	$\sigma_i \gets \sigma_\mathrm{max}, \ i=1,\dots,d$ \tcp*{Initialize standard deviations}
	$r_i \gets \beta \sigma_i, \ i=1,\dots,d$ \tcp*{Initialize range factors}

	$k \gets 1$ \tcp*{Iteration count}

	\While{$k \leq k^\mathrm{max}$}{
	$\mathbf{D} \gets \diag({\{\sigma_i\}}_{i=1}^d), \ \boldsymbol{\Sigma} \coloneqq \mathbf{D}^2$ \tcp*{Covariance matrix}

	$N \gets \alpha d$ \tcp*{Number of samples}

	Generate $N$ independent pseudo-random samples from normal distribution:
	$\mathcal{Y} \gets \mathcal{Y} \cup \{\mathbf{y}_i \sim \mathcal{N}(\ubx^{(k)}, \bSigma)
		\mid x_j^\mathrm{min} - \delta \sigma \leq y_j^i \leq x_j^\mathrm{max} + \delta \sigma, \ j=1,\dots,d\}_{i=1}^{N}$\;

	$\mathcal{Q} \gets \mathcal{Q} \cup \{q_i = f(\mathbf{y}_i)\}_{i=1}^{N} $ \tcp*{Run delamination simulations and compute QoI}

	$\tilde{f} \coloneqq \texttt{NN-Interpolator}(\mathcal{Y}, \mathcal{Q})$ \tcp*{Construct nearest-neighbor (NN) interpolator}
	$\mathcal{R} \coloneqq \left\{ \ubx \in \mathbb{R}^d \mid \max\{x_i^{(k)} - r_i, x_i^\mathrm{min}\} \leq x_i \leq \min\{x_i^{(k)} + r_i,x_i^\mathrm{max} \}, \ i=1,\dots,d \right\}$ \tcp*{RoI}

	$\mathbf{x}^{(k+1)} \gets \argmin_{\ubx \in \mathcal{R}} \tilde{F}(\ubx; \bSigma, \tilde{f})$ \tcp*{Solve optimization sub-problem~\eqref{eq:optimization_subproblem}}

	$\ubx^{*} \gets \ubx^{(k+1)}$ \;

	$s_i^{(k)} \gets \frac{x^{(k+1)}_i - x^{(k)}_i}{r_i},\quad i=1,\dots,d$\ \tikzmark{2nd} \tcp*{Step size}
	\eIf{$k > 1$}{
	$c_i \gets s_i^{(k)} s_i^{(k-1)}$\;
	$\hat{c}_i = \sqrt{\vert c_i \vert} \sign(c_i)$ \tcp*{Oscillation indicator}
	}{
	$\hat{c}_i \gets 1$
	}

	$\gamma_i \gets \frac{1}{2}[\gamma_\mathrm{pan} (1+\hat{c}_i) + \gamma_\mathrm{osc} (1-\hat{c}_i)]$ \quad \tikzmark{right} \tcp*{Contraction parameter}
	$\lambda_i \gets \eta  + \vert s_i^{(k)} \vert (\gamma_i - \eta)$
	\tcp*{Contraction rate}
	$r_i \gets \lambda_i r_i$ \tikzmark{4th} \tcp*{Update range factors}
	$\sigma_i \gets \min\left\{ \max\{\frac{1}{\beta} r_i, \sigma\}, \sigma^{\mathrm{max}} \right\} , \ i=1,\dots,d$ \tcp*{Update standard deviations}

	$k \gets k + 1$\;}%
	\KwRet{$\mathbf{x}^{*}$}
	\begin{tikzpicture}[overlay, remember picture]
		\node[anchor=base] (a) at (pic cs:2nd) {\vphantom{h}}; 
		\node[anchor=base] (b) at (pic cs:4th) {\vphantom{g}}; 
		\draw [decoration={brace,amplitude=0.5em},decorate,ultra thick,gray]
		(a.north -| {pic cs:right}) -- (b.south -| {pic cs:right}) node[midway, right=0.5cm, black] {\emph{Move limit heuristic}~\cite{OnTheRobustneStande2002}};
	\end{tikzpicture}%
\end{algorithm2e}

The optimization algorithm works as follows:
\begin{enumerate}
	\item We start by choosing an initial guess for the design, which is done based on Design of Experiments (DoE)~\cite{EngineeringDesForres2008}.
	      The methods range from simple random sampling to more sophisticated sampling techniques.
	      In this work, we sample $n_0$ uniform random points in the search space $\mathcal{B}$ using Latin Hypercube Sampling (LHS)~\cite{DesignModellingExperimentsFang2005} method, for which we evaluate the QoI.
	      The best design among all serves as the starting point $\ubx^{(1)}$ for the successive surrogate-based optimization procedure.
	      Here, the superscript $(\,\cdot\,)^{(k)}$ denotes the optimization iteration number.
	\item We choose a starting standard deviation parameter $\sigma^\mathrm{max}$ which must be larger than the desired standard deviation $\sigma$, resulting in the starting covariance matrix $\bSigma^{(1)} = {\{\sigma^\mathrm{max}\}}^2 \mathbf{I}$.
	      The covariance matrix $\bSigma$ serves two purposes: (a) in each optimization iteration $k$, it is used to generate random samples from the normal distribution $\mathcal{N}(\ubx^{(k)}, \bSigma^{(k)})$, and (b) it is used to calculate the expectation $\tilde{F}(\ubx; \bSigma^{(k)}, \tilde{f})$, defined in Eq.~\eqref{eq:tilde_F}, for a given design $\ubx$ and approximate response surface $\tilde{f}$.
	      We also define initial range factors $r_i^{(1)} = \beta \sigma^\mathrm{max}, i=1,\dots,d$ for each of the design variables, that allow us to define a RoI centered at $\ubx^{(k)}$, $k=1$:
	      \begin{equation}
		      \mathcal{R}^{(k)} \coloneqq \left\{ \ubx \in \mathbb{R}^d \mid \max\{x_i^{(k)} - r_i^{(k)}, x_i^\mathrm{min}\} < x_i < \min\{x_i^{(k)} + r_i^{(k)},x_i^\mathrm{max} \}, \ i=1,\dots,d \right\},
		      \label{eq:region_of_interest}
	      \end{equation}
	      where $\beta$ is a hyperparameter, ideally selected between \num{1} and \num{3}.
	\item In each optimization iteration, we generate $n=\alpha d$ samples randomly from the normal distribution with mean $\ubx^{(k)}$ and covariance $\bSigma^{(k)}$, and evaluate the structural response $f$ for each of the new samples.
	      Here, $\alpha$ is the ``over-sampling'' factor, chosen between \num{2} and \num{10}.
	      All the available data combined with the data from the previous optimization iterations is then used to construct a nearest-neighbor interpolant $\tilde{f}$.
	\item Given a region of interest $\mathcal{R}^{(k)}$, the covariance matrix $\bSigma^{(k)}$ and the response function approximation $\tilde{f}$ for an optimization iteration $k$, we solve the following optimization sub-problem to obtain the next design point $\ubx^{(k+1)}$:
	      \begin{align}
		      \min_{\mathbf{x} \in \mathcal{R}^{(k)}}\ \tilde{F}(\ubx; \bSigma^{(k)}, \tilde{f}).
		      \label{eq:optimization_subproblem}
	      \end{align}
	      Since, $\tilde{F}$ is a smooth function of $\ubx$, the above optimization sub-problem can in principle be solved using a gradient-based solver, for example, based on trust-region methods~\cite{TrustRegionMeConn2000}.
	      The gradient and the Hessian of $\tilde{F}$ with respect to $\ubx$ are estimated using Monte Carlo evaluations of $\tilde{f}$, as explained in the next section.
	      Due to the limitations on the accuracy of estimators utilizing Monte Carlo methods, we perform only few iterations of the optimization solver to obtain an approximate solution.
	\item For updating the RoI and the covariance matrix, we use the \emph{move limit heuristic} (see~\cite{OnTheRobustneStande2002} for details), which updates the range factors based on the distance between the new design $\ubx^{(k+1)}$ and the old design $\ubx^{(k)}$.
	      Since, we approximate the response surface by nearest-neighbor interpolation, this can result in severe oscillations between the successive designs.
	      To prevent this undesired behavior, the heuristic uses oscillation indicators
	      \begin{align}
		      \hat{c}_i = \sqrt{\vert c_i \vert} \sign(c_i), \quad \text{with}
		      \quad c_i = s_i^{(k)} s_i^{(k-1)}, \quad i=1,\dots,d,
	      \end{align}
	      where $s_i^{(k)}$ is the design step size, normalized w.r.t.\ the range factor $r_i^{(k)}$, expressed as
	      \begin{align}
		      s_i^{(k)} = \frac{x^{(k+1)}_i - x^{(k)}_i}{r_i^{(k)}}.
	      \end{align}
	      As part of the move-limit heuristic, the contraction parameters are then defined by
	      \begin{align}
		      \gamma_i = \frac{1}{2}[\gamma_\mathrm{pan} (1+\hat{c}_i) + \gamma_\mathrm{osc} (1-\hat{c}_i)], \quad i=1,\dots,d,
	      \end{align}
	      where $\gamma_\mathrm{pan}$ and $\gamma_\mathrm{osc}$ are the hyperparameters, accounting for the directional movement and oscillations in the design, respectively.
	      With zoom parameter $\eta$, the contraction rates are defined as
	      \begin{align}
		      \lambda_i = \eta  + \vert s_i^{(k)} \vert (\gamma_i - \eta),\quad i=1,\dots,d,
	      \end{align}
	      which results to updates of the range factors:
	      \begin{align}
		      r_i^{(k+1)} = \lambda_i r_i^{(k)},\quad i=1,\dots,d.
	      \end{align}
	      The updated range factors are then used to define the RoI for the next optimization iteration.
	      We also update the covariance to a diagonal matrix:
	      \begin{align}
		      \bSigma^{(k+1)} = \diag(\{\sigma_1^2,\dots,\sigma_d^2\}),
		      \label{eq:update_covariance_matrix}
	      \end{align}
	      where $\sigma_i, i=1,\dots,d$ are the internal parameters of the algorithm, acting as intermediate standard deviations of the design parameters, given by
	      \begin{equation}
		      \sigma_i = \min\left\{ \max\{\frac{1}{\beta} r_i^{(k+1)}, \sigma\}, \sigma^{\mathrm{max}} \right\},\quad \ i=1,\dots,d.
	      \end{equation}
	\item We choose the stopping criterion based on the computational budget, which defines the maximum number of optimization iterations $k^\mathrm{max} \in \mathbb{N}$.
	      Also, the maximum number of iterations should be large enough such that $\sigma_i,\ i=1,\dots,d$ are reduced to the desired standard deviation $\sigma$ at the end of an optimization run.
\end{enumerate}

\subsubsection{Search gradients for optimization sub-problem}%
\label{sec:search_gradients}
A good property of using the multivariate normal distribution for the stochastic optimization problem is that the probability density function is at least twice differentiable, allowing computation of its first and second derivatives with respect to the mean.
The gradient and the Hessian of $p$ with respect to the mean $\ubx$ are given by
\begin{align}
  \nabla_{\ubx} p(\by; \ubx, \bSigma) &= p(\by; \ubx, \bSigma)\bSigma^{-1} [\by-\ubx], \\
  \nabla^2_{\ubx} p(\by; \ubx, \bSigma) &= p(\by; \ubx, \bSigma) \left[-\bSigma^{-1} + \bSigma^{-1}[\by-\ubx][\by-\ubx]^\top \bSigma^{-1}\right],
\end{align}
respectively.
Using the above information, we can define the gradient and the Hessian of $\tilde{F}$ with respect to the design $\ubx$, as
\begin{align}
	\nabla_\ubx \tilde{F}(\ubx; \bSigma, \tilde{f}) & = \int_{\mathbb{R}^d} \tilde{f} (\by) \nabla_\ubx p(\by; \ubx, \bSigma)\, \dd{\by}             
	                                                = \int_{\mathbb{R}^d} \tilde{f}(\by) \bSigma^{-1} [\by - \ubx] p(\by; \ubx, \bSigma)\, \dd{\by} \nonumber \\
	                                                & = \mathbb{E}_{\ubx, \bSigma} \left[\tilde{f}(\bY) \bSigma^{-1} [\bY - \ubx]\right], \\
	\nabla^2_\ubx \tilde{F}(\ubx; \bSigma, \tilde{f}) & = \int_{\mathbb{R}^d} \tilde{f} (\by) \nabla^2_\ubx p(\by; \ubx, \bSigma)\, \dd{\by}
	                                                   = \int_{\mathbb{R}^d} \tilde{f}(\mathbf{y}) \left[-\boldsymbol{\Sigma}^{-1} + \boldsymbol{\Sigma}^{-1} [\mathbf{y} - \mathbf{x}][\mathbf{y} - \mathbf{x}]^\top \boldsymbol{\Sigma}^{-1}\right] p(\by;\ubx, \bSigma) \, \dd{\by}                   \nonumber \\
	                                                  & = \mathbb{E}_{\ubx, \bSigma} \left[ \tilde{f}(\mathbf{Y}) \left[-\boldsymbol{\Sigma}^{-1} + \boldsymbol{\Sigma}^{-1} [\mathbf{Y} - \mathbf{x}][\mathbf{Y} - \mathbf{x}]^\top \boldsymbol{\Sigma}^{-1}\right] \right],
\end{align}
respectively.
The $d$-dimensional expectations in the above derivatives can be evaluated approximately or with arbitrary precision using Monte Carlo methods.
Given $M$ normally distributed random samples $\bz_i \sim \mathcal{N}(\ubx, \bSigma),\ i=1,\dots,M$, the sample average approximations of the objective and its derivatives are given by
\begin{subequations}
	\begin{align}
		\tilde{F}(\ubx; \bSigma, \tilde{f})               & \approx   \frac{1}{M} \sum_{i=1}^M \tilde{f}(\bz_i),     \label{eq:sample_average_objective}                     \\
		\nabla_\ubx \tilde{F}(\ubx; \bSigma, \tilde{f})   & \approx \frac{1}{M} \sum_{i=1}^{M} \tilde{f}(\mathbf{z}_i) \boldsymbol{\Sigma}^{-1} [\mathbf{z}_i - \mathbf{x}],
		\label{eq:sample_average_objective_gradient}                                                                                                                         \\
		\nabla^2_\ubx \tilde{F}(\ubx; \bSigma, \tilde{f}) &
		\approx \frac{1}{M} \sum_{i=1}^{M} \tilde{f}(\mathbf{z}_i) \left[-\boldsymbol{\Sigma}^{-1} + \boldsymbol{\Sigma}^{-1} [\mathbf{z}_i - \mathbf{x}][\mathbf{z}_i - \mathbf{x}]^\top \boldsymbol{\Sigma}^{-1}\right].
		\label{eq:sample_average_objective_hessian}
	\end{align}
	\label{eq:sample_average}
\end{subequations}
In practice, it has been observed that using quasi-random samples, e.g., generated from low-discrepancy Sobol sequences~\cite{ConstructingSoJoeS2008}, instead of the usual pseudo-random points leads to a better convergence rate in the Monte Carlo estimation of the above expectations.
For the numerical problems discussed in the next section, we used $M = 2^{16} = \num{65536}$ samples.

\section{Numerical experiments}%
\label{sec:numerical_experiments}

In this section, we present an example of design optimization for improving delamination resistance of a heterogeneous double-cantilever beam using the algorithm described in the previous section.
The double-cantilever beam setup is the most commonly applied test to measure mode-\rom{1} interlaminar fracture toughness~\cite{TestingTheTouLaffan2012}.
Here, we choose the design parameters to be the shapes of the elliptical, comparatively hard inclusions embedded in the beam (matrix) with predefined constitutive properties of the bulk as well as the material interfaces.
The objective is to minimize the expected QoI defined by~\eqref{eq:mechanical_work} over a chosen standard deviation parameter $\sigma$, leading to a stochastic optimization problem specified by Definition~\ref{def:optimization_problem}.

\begin{figure}[th]
	\centering
	\begin{tikzpicture}[scale=0.9, text=black, draw=black]
	\draw[fill=black!2] (-5.5, -1.3) rectangle (9.5, 1.3);
	\draw[blue, dashed, thick] (-5.5, 0) -- (0, 0);
	\draw[red, thick] (0, 0) -- (9.5, 0);

	\draw[-latex] (-0.6, 0) to[out=-90, in=0] (-1.2, -0.5) node[anchor=east] {\emph{pre-crack}};

	\draw[fill=green!20, draw=red, thick, rotate around={45:(0.8, 0)}] (0.8, 0) ellipse (0.25 and 0.15);
	\draw[fill=green!20, draw=red, thick, rotate around={10:(0.4, 0.4)}] (0.4, 0.4) circle (0.25 and 0.18);
	\draw[fill=green!20, draw=red, thick, rotate around={60:(0.4, -0.4)}] (0.4, -0.4) circle (0.23 and 0.15);
	\draw[fill=green!20, draw=red, thick, rotate around={47:(1.6, 0)}] (1.6, 0) circle (0.25 and 0.25);
	\draw[fill=green!20, draw=red, thick, rotate around={90:(1.2, 0.4)}] (1.2, 0.4) circle (0.2 and 0.18);
	\draw[fill=green!20, draw=red, thick, rotate around={65:(1.2, -0.4)}] (1.2, -0.4) circle (0.25 and 0.2);
	\draw[fill=green!20, draw=red, thick, rotate around={5:(2.4, 0)}] (2.4, -0.) circle (0.25 and 0.16);
	\draw[fill=green!20, draw=red, thick, rotate around={1:(2, 0.4)}] (2., 0.4) circle (0.25 and 0.23);
	\draw[fill=green!20, draw=red, thick, rotate around={25:(2, -0.4)}] (2., -0.4) circle (0.25 and 0.21);
	\draw[fill=green!20, draw=red, thick, rotate around={18:(3.2, 0)}] (3.2, -0.) circle (0.25 and 0.18);
	\draw[fill=green!20, draw=red, thick, rotate around={63:(2.8, 0.4)}] (2.8, 0.4) circle (0.25 and 0.15);
	\draw[fill=green!20, draw=red, thick, rotate around={110:(2.8, -0.4)}] (2.8, -0.4) circle (0.15 and 0.15);

	\def\s{0.05}
	\foreach \x in {0.8, 1.6, 2.4, 3.2} {

			\draw[ultra thin] (\x-\s, 0) -- (\x+\s, 0);
			\draw[ultra thin] (\x, -\s) -- (\x, \s);

			\draw[ultra thin] (\x-0.4-\s, 0.4) -- (\x-0.4+\s, 0.4);
			\draw[ultra thin] (\x-0.4, -\s+0.4) -- (\x-0.4, \s+0.4);

			\draw[ultra thin] (\x-0.4-\s, -0.4) -- (\x-0.4+\s, -0.4);
			\draw[ultra thin] (\x-0.4, -\s-0.4) -- (\x-0.4, \s-0.4);

		}

	\draw[dashed, ultra thin] (-0.2, -0.8) rectangle (3.8, 0.8);
	\draw[-latex] (2.3, -0.8) to[out=-90, in=180] (4, -1.7) node[anchor=west] {\emph{heterogeneous region}};

	\draw[-latex, brown, ultra thick] (-6, -0.1) -- (-6, -1.2) node[midway, left, black] {$-u_2$};
	\draw[-latex, brown, ultra thick] (-6, 0.1) -- (-6, 1.2) node[midway, left, black] {$u_2$};

	\foreach \y in {1.2,1,...,-1.2}
	\draw (9.5, \y) -- (9.8, \y-0.2);

	\dimline[extension start length=0.2cm, extension end length=0.2cm] {(-5.5, 1.6)}{(9.5, 1.6)}{$15$};
	\dimline[extension start length=0.4cm, extension end length=0.4cm] {(10.1, 1.3)}{(10.1, 0)}{$1.3$};
	\dimline[extension start length=0.4cm, extension end length=0.4cm] {(10.1, 0)}{(10.1, -1.3)}{$1.3$};
	\dimline[label style={fill=black!2}, extension start length=0cm, extension end length=0.3cm] {(-5.5, 0.3)}{(0, 0.3)}{$5.5$};
	\dimline[label style={above, fill=none}, extension start length=1.1cm, extension end length=1.1cm] {(4, 0.4)}{(4, -0.4)}{$0.8$};
	\dimline[label style={fill=none, above}, extension start length=0.3cm, extension end length=0.3cm] {(2., 0.7)}{(2.8, 0.7)}{$0.8$};
	\dimline[label style={fill=none, above}, extension start length=0.7cm, extension end length=0.7cm] {(0., 0.7)}{(0.8, 0.7)}{$0.8$};

\end{tikzpicture}
	\caption{
		Schematic showing the problem setup for shape optimization, where the design parameters are the radii and the orientations of the elliptical inclusions.
		The initial crack is shown in blue, while the interfaces in red represent the potential crack paths.
	}
	\label{fig:problem_setup}
\end{figure}
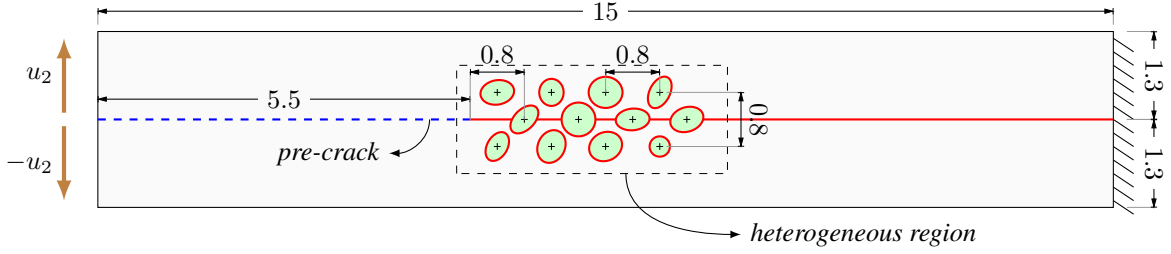

\begin{figure}
	\centering
	\begin{minipage}{0.5\textwidth}
		\centering
		\small{
	\def\svgwidth{1\columnwidth}
	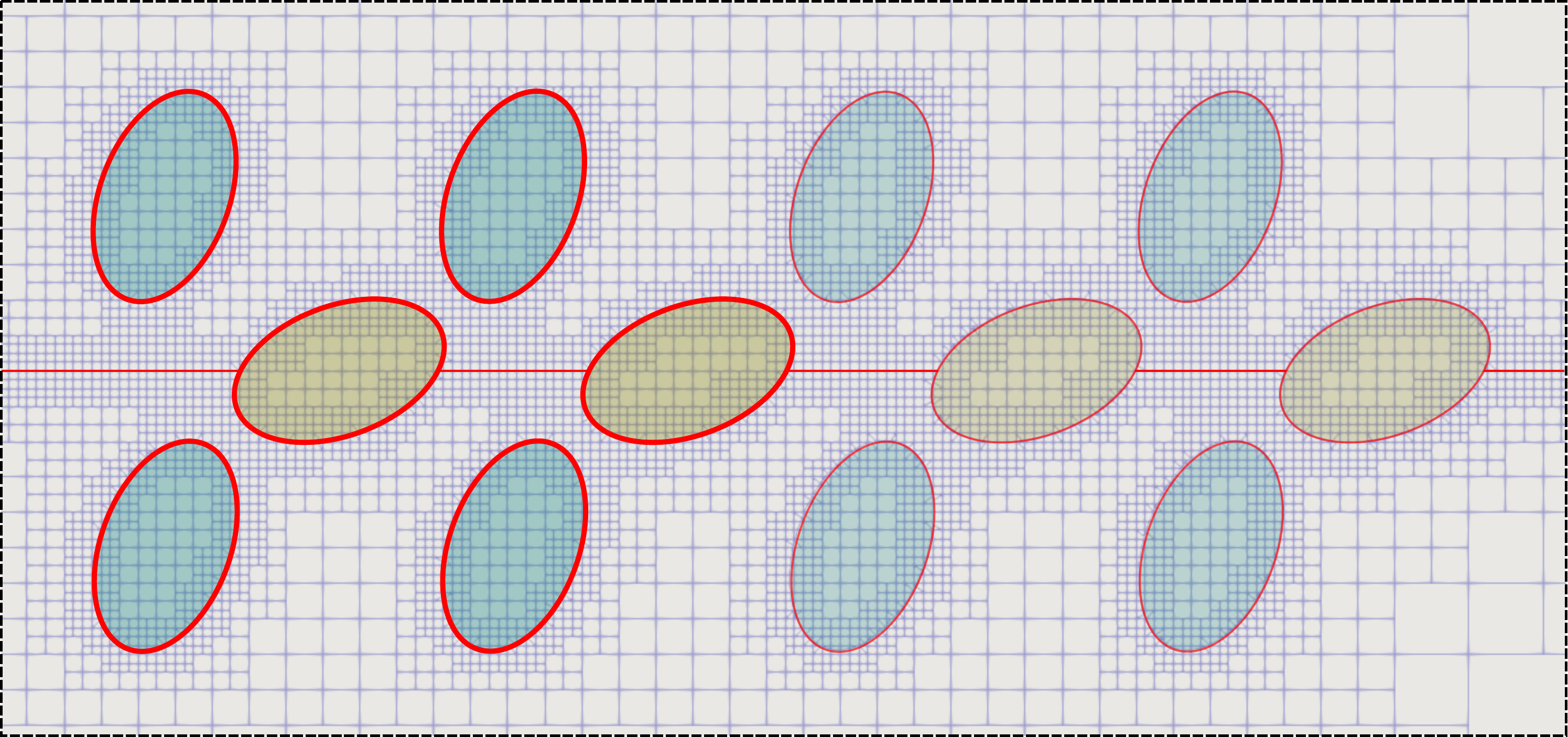
}
	\end{minipage}\hfill
	\begin{minipage}{0.48\textwidth}
		\centering
		\begin{tabular}{ c | c }
			\toprule
			$d$  & Design parameterization                                                 \\
			\midrule
			\midrule
			$2$  &
			$x_1=\theta_1 = \theta_4, \quad x_2=\theta_2 = \theta_3 = \theta_5 = \theta_6$ \\
			\hline
			$6$  &
			\begin{tabular}{c c c}
				$x_1=\theta_1=\theta_4$, & $x_2=\theta_2=\theta_5$, & $x_3=\theta_3=\theta_6$ \\
				$x_4=r_1=r_4$,           & $x_5=r_2=r_5$,           & $x_6=r_3=r_6$
			\end{tabular}  \\
			\hline
			$12$ &
			\begin{tabular}{c c}
				$x_i=\theta_i,\ i=1,\dots,6$ \\
				$x_{i+6}=r_i,\ i=1,\dots,6$
			\end{tabular}
			\\
			\bottomrule
		\end{tabular}
	\end{minipage}
	\captionof{figure}{
	(Left) Structural setup with $12$ shape parameters $\{\theta_i, r_i\}_{i=1}^6$.
	The radii along the minor axes of the elliptical inclusions are fixed at \num{0.15}, while the radii $r_i, i=1,\dots,d$ along the major axes range between \num{0.15} and \num{0.25}.
	The orientations $\theta_i, i=1,\dots,6$ can take any values in $\mathbb{R}$.
	(Right) Parameterization for problem dimensions $d=2, 6$ and $12$ with the design variables $x_i, i=1,\dots,d$.
	}
	\label{fig:design_parameters}
\end{figure}

The geometric setup and boundary conditions are shown in Fig.~\ref{fig:problem_setup}.
With the right boundary fixed, the left boundary is subjected to a monotonically increasing displacement $\hat{u}_2(\,\cdot\,,t) = \pm t, \ t \in [0,\ T=0.1]$ in the vertical direction, while it is free to move horizontally.
The loading is applied using $N_T = \num{100}$ equal time steps, resulting in a load/time step size of $\Delta t = \num{0.001}$.
We assume isotropic material for the two beams, upper and lower, with Young's modulus $E=\num{200}$ and Poisson's ratio $\nu=\num{0.3}$ (or in terms of Lam\'e parameters: $\lambda=\num{115.38}$ and $\mu=76.92$).
Similarly, for the inclusions, we assume $E=\num{400}$ and $\nu=\num{0.3}$ (or $\lambda=\num{230.77}, \mu=\num{153.85}$).
The matrix-matrix and matrix-inclusion interfaces are modeled using the exponential cohesive law, as described by Eq.~\eqref{eq:needleman_fracture_energy}, with parameters: $r=\num{0}$, $\delta_\tn^*=\delta_\ts^*=\num{0.0001}$, $\phi_\tn=\num{2.718e-05}$, and $\phi_\ts=\num{1.166e-5}$.
The spatial domain is discretized with Lagrange finite elements of polynomial degree \num{1}.
The resulting finite element problem was implemented using  an open-source finite element library, \texttt{deal.II}~\cite{TheDealIiFinArndt2021}.

In the following optimization studies, the algorithmic parameters in Algorithm~\ref{alg:optimization_algorithm} are chosen as $\alpha = 3$, $n_0 = 3 d$, $\beta = 2$, $\delta=3$, $\gamma^\mathrm{pan} = 1.2$, $\gamma^\mathrm{osc} = 0.8$, $\eta = 0.8$, and $k^{\mathrm{max}} = 100$.

\setlength{\fboxsep}{0pt}
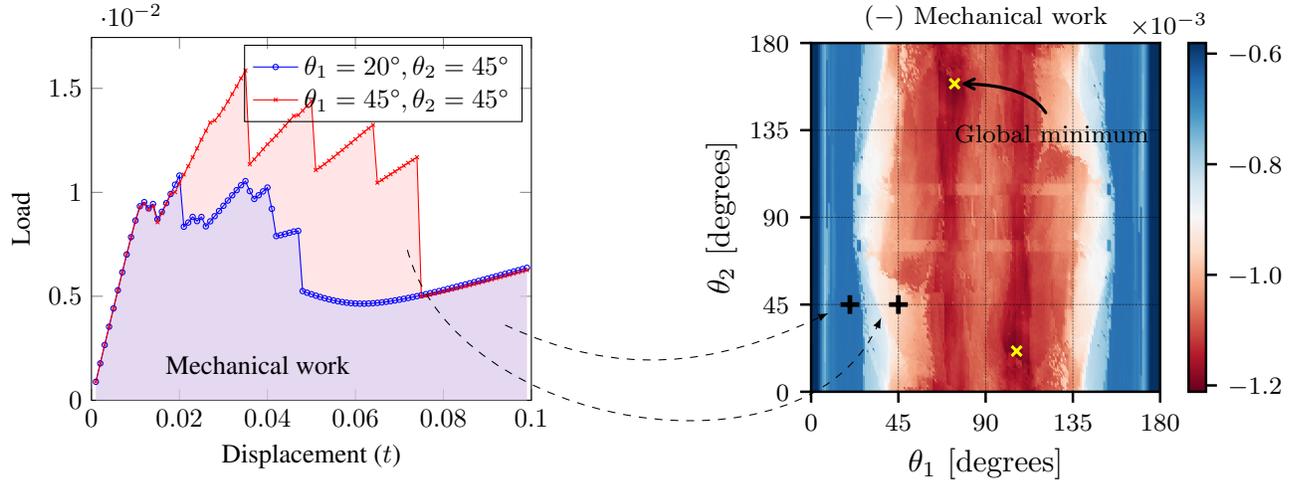
\begin{figure}
	\begin{center}
		\begin{tikzpicture}[text=black, draw=black]
			\begin{axis} [
					width=0.45\textwidth,
					xlabel={Displacement ($t$)},
					ylabel={Load},
					xlabel near ticks,
					ylabel near ticks,
					xmin = 0, xmax=0.1,
					ymin = 0,
					xticklabel style={
							/pgf/number format/precision=3,
							/pgf/number format/fixed,
						},
					legend style = {fill=none,},
				]
				\path [name path=xaxis] (\pgfkeysvalueof{/pgfplots/xmin},0) -- (\pgfkeysvalueof{/pgfplots/xmax},0);
				\addplot[mark=o, draw=blue, mark size=1pt, name path=f2] table[header=false, x index=1, y index=7] {data/solution_statistics_2d_20_45.dat} node[pos=0.4, below=2cm] {Mechanical work};
				\addlegendentry{$\theta_1=\SI{20}{\degree}, \theta_2=\SI{45}{\degree}$}
				\addplot[mark=x, draw=red, mark size=1pt, name path=f1] table[header=false, x index=1, y index=7] {data/solution_statistics_2d_45_45.dat};
				\addlegendentry{$\theta_1=\SI{45}{\degree}, \theta_2=\SI{45}{\degree}$}

				\addplot[red!20, opacity=0.5] fill between[of=f1 and xaxis, soft clip={domain=0.001:0.099}];
				\addplot[blue!20, opacity=0.5] fill between[of=f2 and xaxis, soft clip={domain=0.001:0.099}];

			\end{axis}
			\node (landscape) at (12, 2.1) {\includegraphics[width=0.5\textwidth]{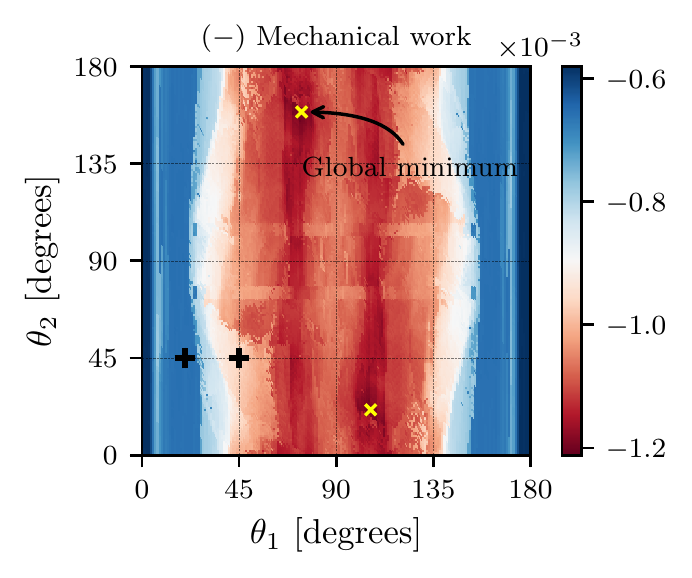}};
			\draw[-latex, thin, dashed] (5.5, 1) to[in=205, out=-25] (9.8, 1.16);
			\draw[-latex, thin, dashed] (4.2, 2) to[in=245, out=-75] (10.5, 1.1);
		\end{tikzpicture}
	\end{center}
	\caption{
		(Left) Load--displacement curves for $\theta_1=\SI{20}{\degree}, \theta_2=\SI{45}{\degree}$ (blue) and $\theta_1=\theta_2=\SI{45}{\degree}$ (red), with the mechanical work indicated by the shaded region.
		(Right) Global landscape of the QoI ($-$ mechanical work) for $\theta_1, \theta_2 \in [\SI{0}{\degree}, \SI{180}{\degree}]$, computed on a regular grid of $180 \times 180$ points.
		(For more information about the color references in this figure, the reader is referred to the digital version of this article.)
	}
	\label{fig:qoi_landscape_2d}
\end{figure}

\begin{figure}
  \begin{center}
    \includegraphics[width=0.5\textwidth]{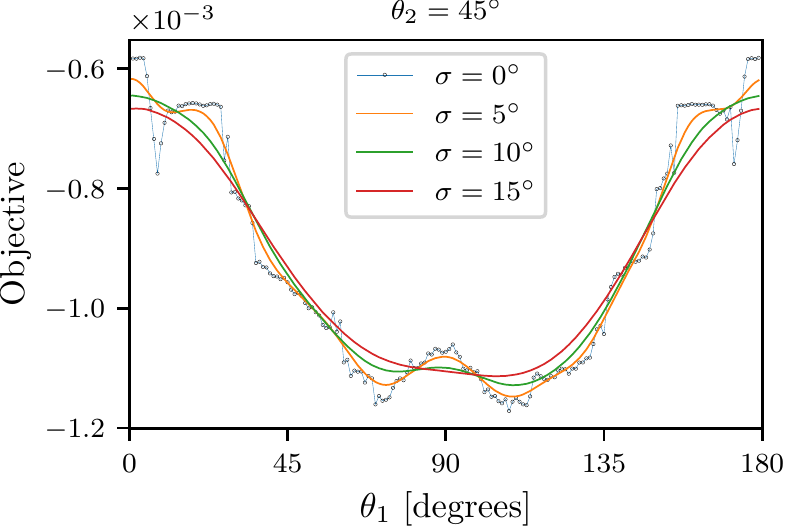}
  \end{center}
  \caption{Objective function landscape in one dimension ($\theta_2 = \SI{45}{\degree}$), with different values of smoothing parameter $\sigma$.
    The non-smoothness of the objective for $\sigma = 0^\circ$ is attributed to finite-element remeshing, presence of structural instabilities and bifurcations.
		(For more information about the color references in this figure, the reader is referred to the digital version of this article.)
  }
  \label{fig:objective_landscape_1d}
\end{figure}

\begin{figure}
	\centering
	\begin{subfigure}{\textwidth}
		\centering
		\includegraphics[width=0.4\textwidth]
		{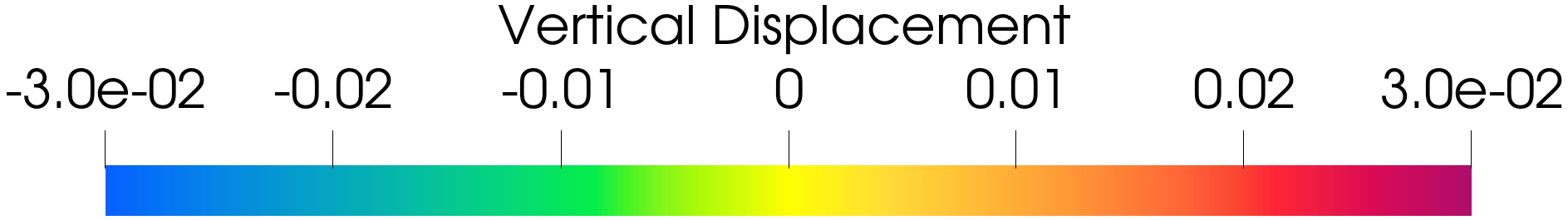}
    \vspace{1em}
	\end{subfigure}
	\begin{subfigure}[t]{0.27\textwidth}
		\centering
		\makebox[0pt][r]{\makebox[25pt]{\raisebox{35pt}{\rotatebox[origin=c]{90}{\begin{tabular}{c} $\theta_1=\SI{20}{\degree}$ \\ $\theta_2=\SI{45}{\degree}$ \end{tabular}}}}}%
		\fbox{\includegraphics[width=\textwidth]
    {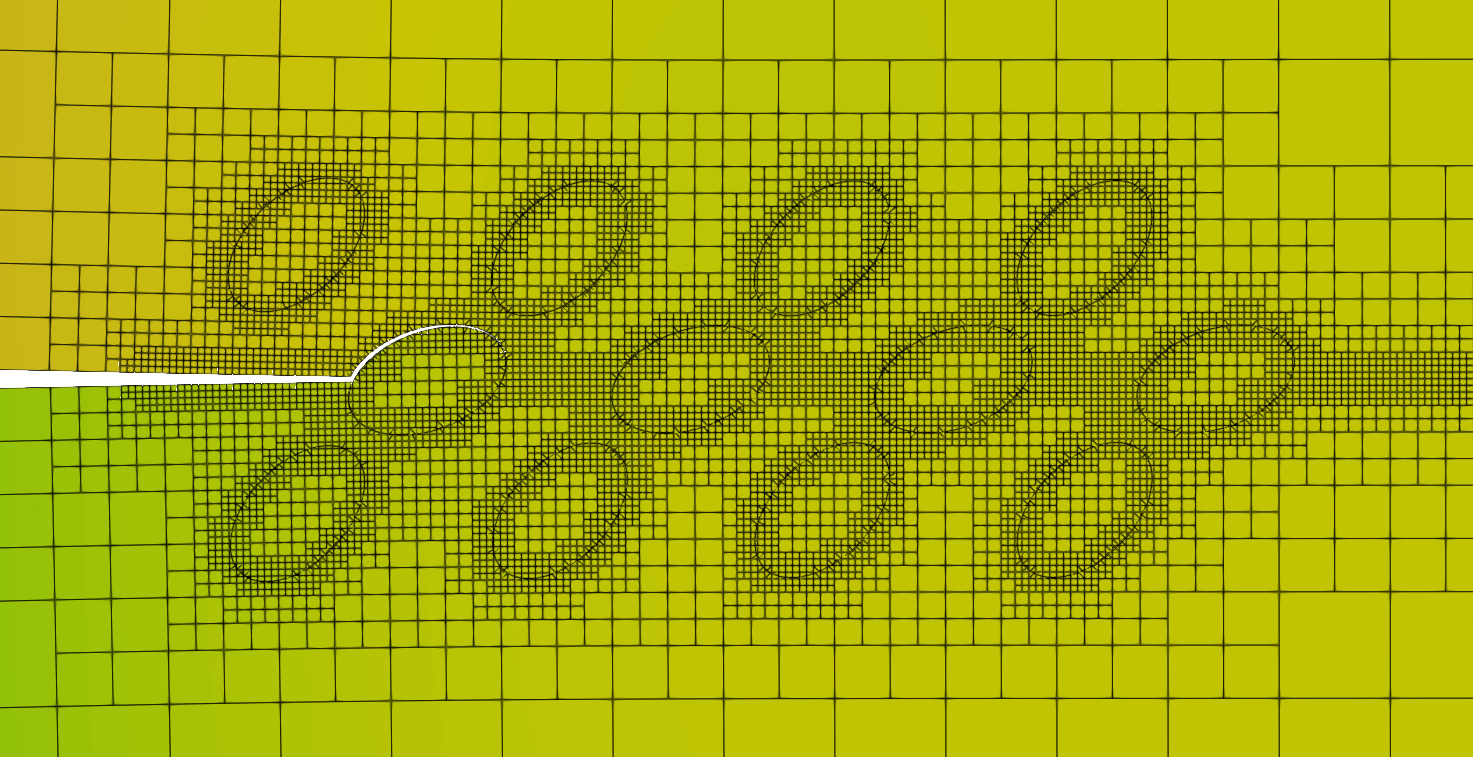}} \vskip 0.1cm
		\makebox[0pt][r]{\makebox[25pt]{\raisebox{35pt}{\rotatebox[origin=c]{90}{\begin{tabular}{c} $\theta_1=\SI{45}{\degree}$ \\ $\theta_2=\SI{45}{\degree}$ \end{tabular}}}}}%
		\fbox{\includegraphics[width=\textwidth]
			{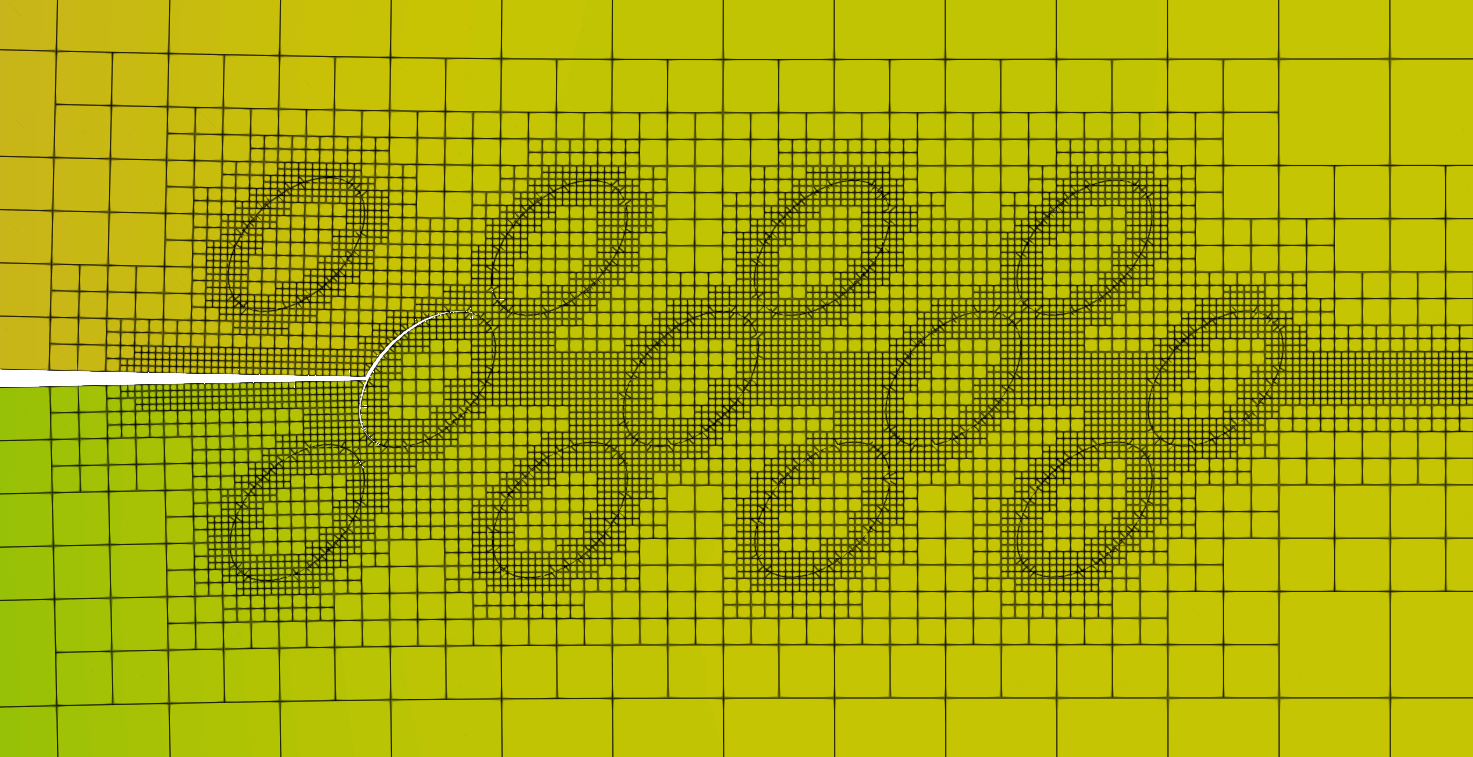}}
      \caption*{\normalsize{$t=0.020$}}
	\end{subfigure}
	\hspace{1em}
	\begin{subfigure}[t]{0.27\textwidth}
		\centering
		\fbox{\includegraphics[width=\textwidth]
			{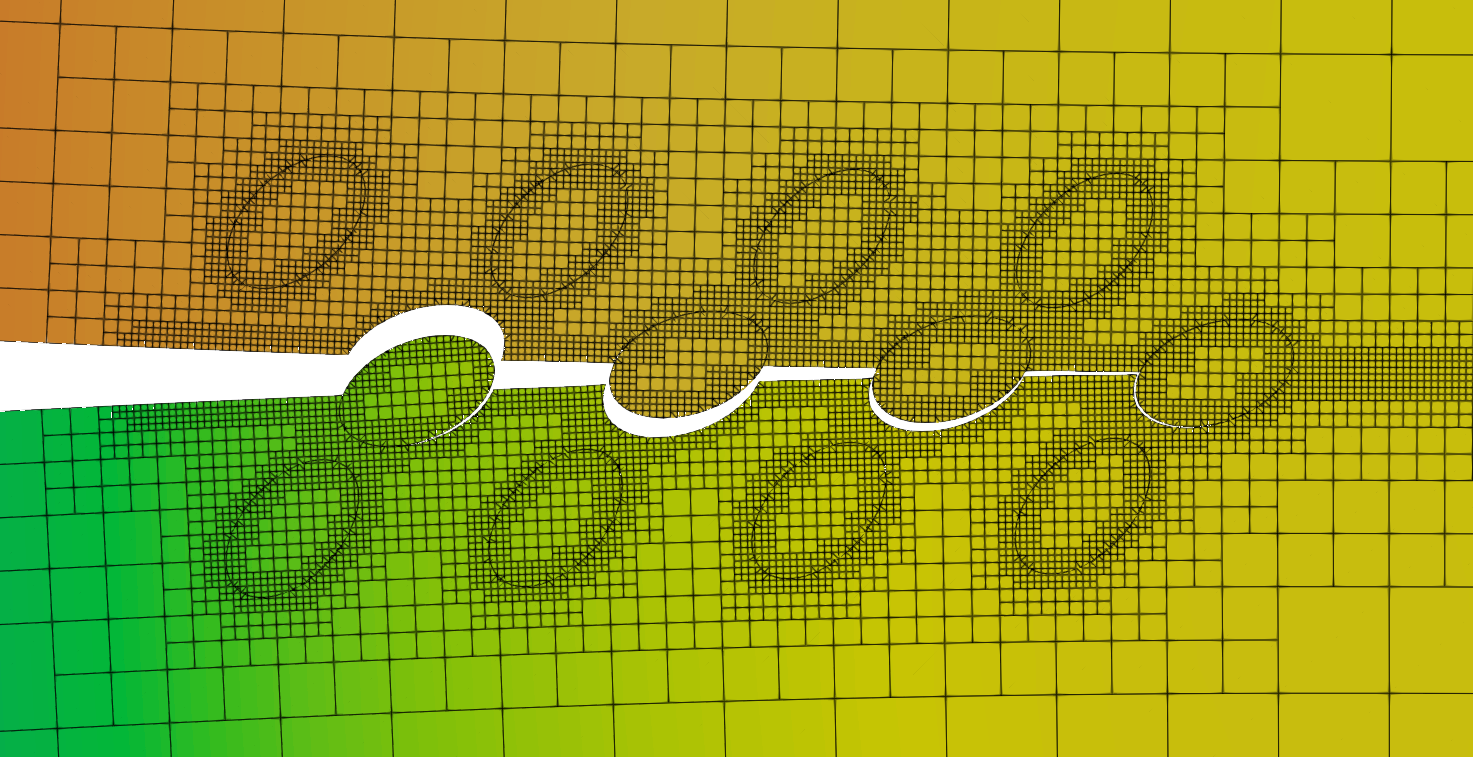}} \vskip 0.1cm
		\fbox{\includegraphics[width=\textwidth]
			{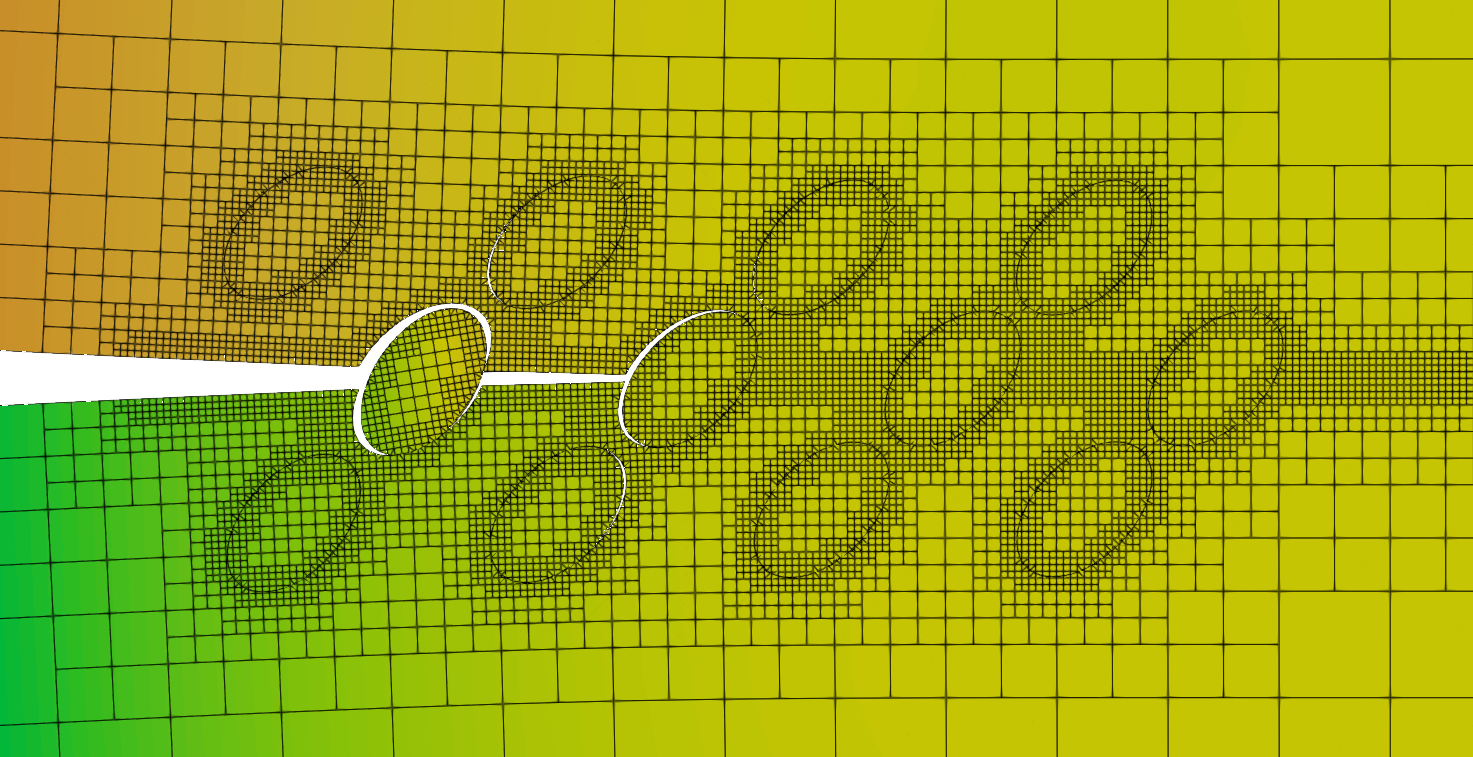}}
      \caption*{\normalsize{$t=0.040$}}
	\end{subfigure}
	\hspace{1em}
	\begin{subfigure}[t]{0.27\textwidth}
		\centering
		\fbox{\includegraphics[width=\textwidth]
			{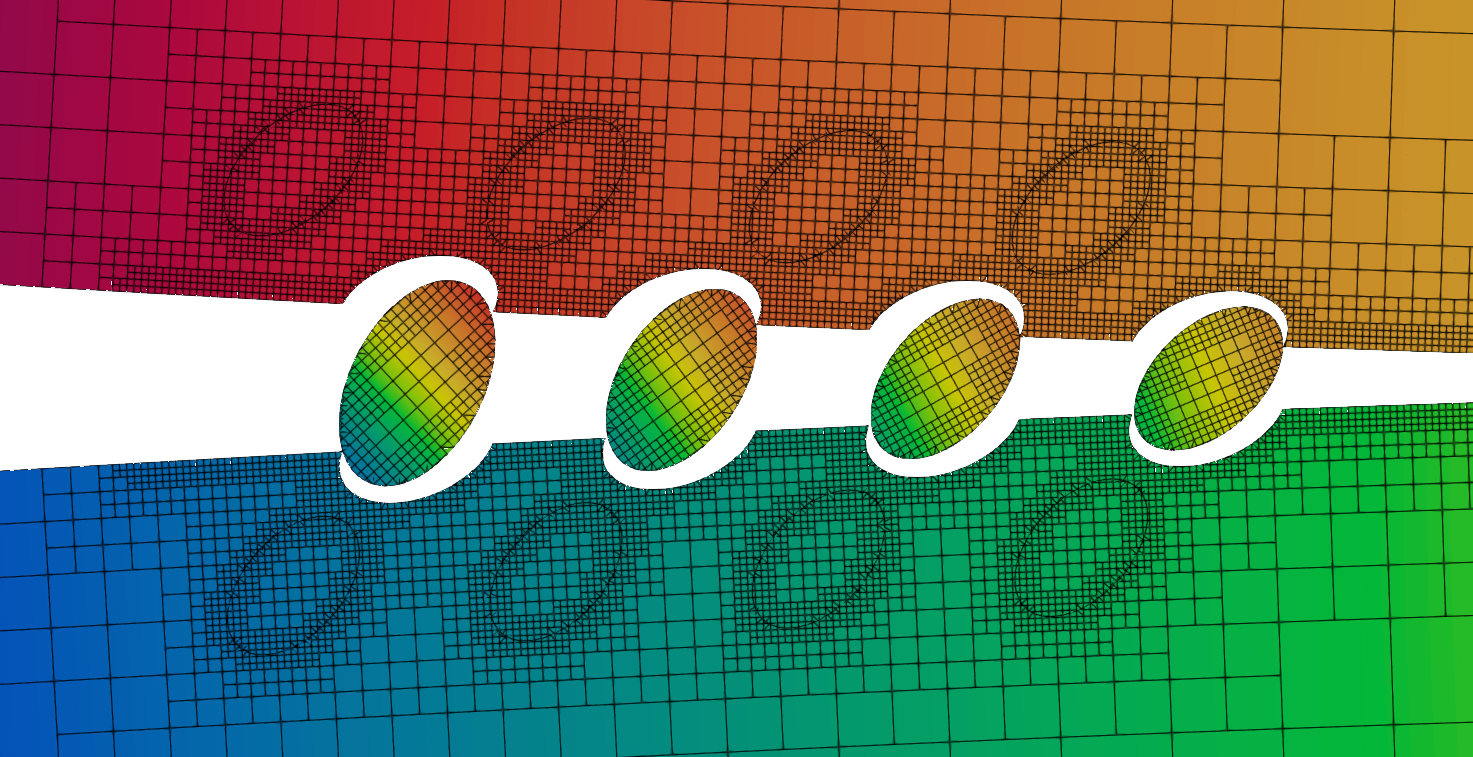}} \vskip 0.1cm
		\fbox{\includegraphics[width=\textwidth]
			{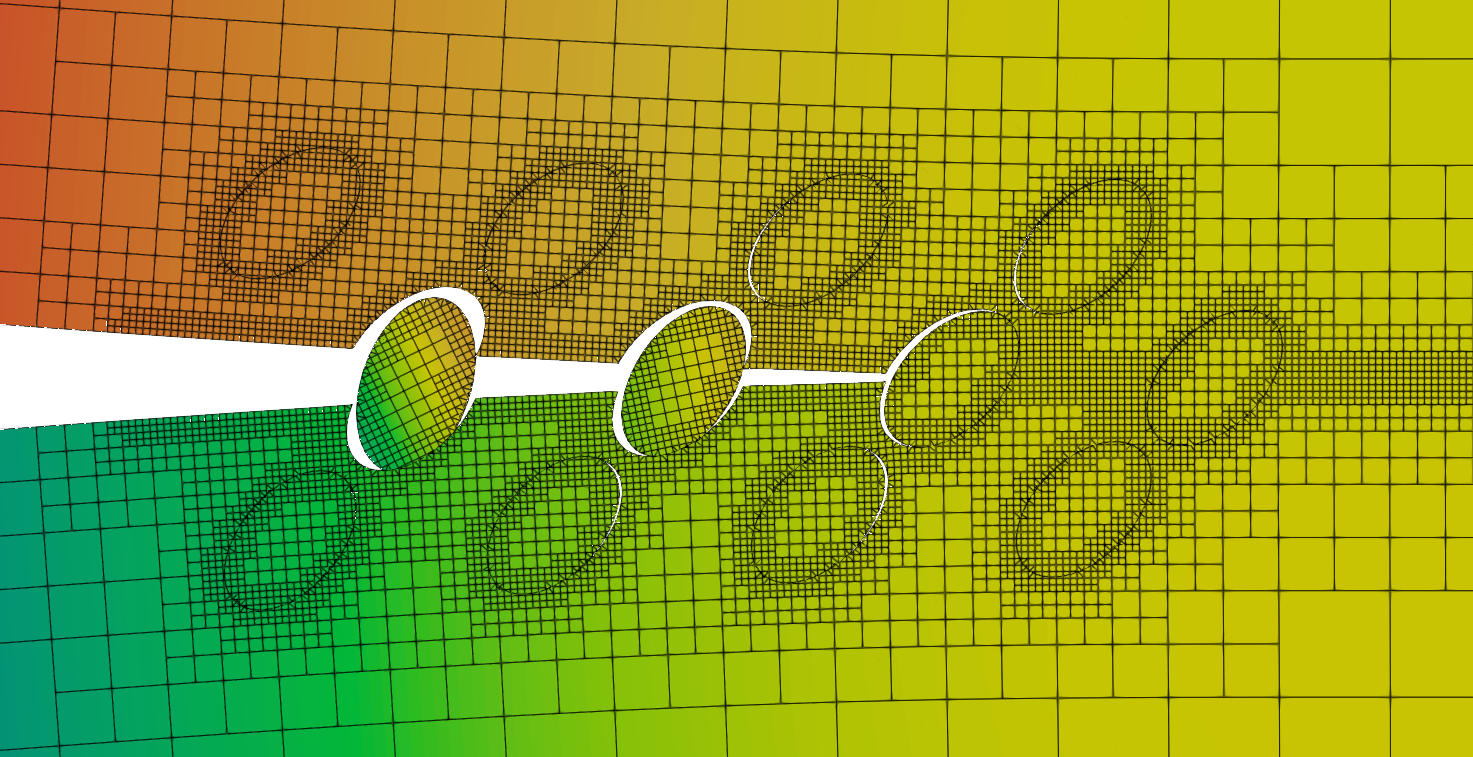}}
      \caption*{\normalsize{$t=0.060$}}
	\end{subfigure}
	\caption{Evolution of cracks for $\theta_1=\SI{20}{\degree}, \theta_2=\SI{45}{\degree}$ (top) and $\theta_1=\SI{45}{\degree}, \theta_2=\SI{45}{\degree}$ (bottom).
		The deformations are scaled by a factor of $10$ for better visualization of the crack openings.
		(For more information about the color references in this figure, the reader is referred to the digital version of this article.)
	}
	\label{fig:visualization_2d_crack}
\end{figure}

\subsection{Optimization with \num{2} DoFs using point cloud data}
\begin{figure}
	\begin{center}
		\includegraphics[width=\textwidth]{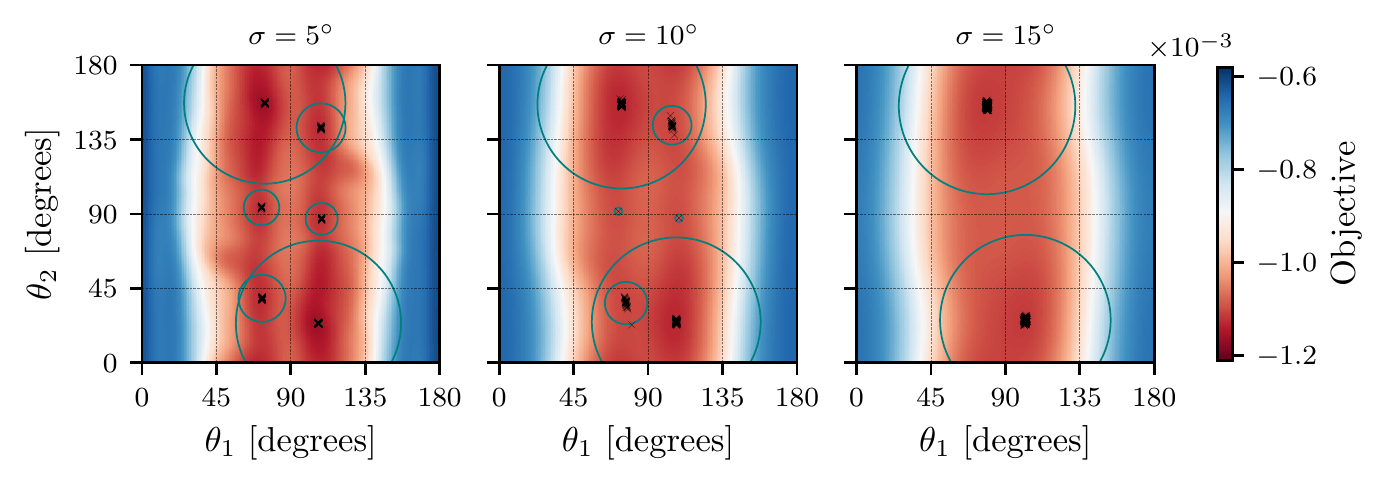}
	\end{center}
	\vspace*{-5mm}
	\caption{
		Global landscapes of the objective function for $\sigma=\SI{5}{\degree}, \SI{10}{\degree}$ and $\SI{15}{\degree}$.
		The scattered points correspond to the final designs obtained by \num{1000} optimization runs, and the size (area) of the enclosing circles is proportional to the number of design points clustered together.
		(For more information about the color references in this figure, the reader is referred to the digital version of this article.)
	}
	\label{fig:objective_landscape_2d_smoothed}
\end{figure}

To get an idea of the behavior of the global landscape of the chosen QoI and the effect of the smoothing parameter $\sigma$ on the optimization procedure over many runs, we first consider a simple case with \num{2} design DoFs ($d=2$).
The design parameters are the angles (in degrees) $\theta_1$ and $\theta_2$, which define the orientations of the two groups of inclusions: those that are located along the central interface of the double cantilever beam, and the rest (see Fig.~\ref{fig:design_parameters}).
The major and minor radii of the elliptical inclusions are set to \num{0.25} and \num{0.15}, respectively
A close approximation to the response surface $f$ is given by a piecewise-constant function
$g(\theta_1, \theta_2) \coloneqq f(\{\lfloor \theta_1 \rfloor , \lfloor \theta_2 \rfloor\})$
generated by evaluating QoI on a regular grid of $180 \times 180$ points for $\theta_1,\theta_2 \in [\SI{0}{\degree}, \SI{180}{\degree}]$.
Fig.~\ref{fig:qoi_landscape_2d} shows the landscape of the QoI, while Fig.~\ref{fig:objective_landscape_1d} depicts a one-dimensional cross-section showing noisy and discontinuous characteristics of the QoI.
The numerical noise stems from the finite-element remeshing for each of the configurations of the inclusions.
Fig.~\ref{fig:visualization_2d_crack} depicts the evolution of cracks for two close configurations of the structural design with a large variation in the QoI.

\begin{figure}[t]
	\begin{center}
		\includegraphics[width=\textwidth]{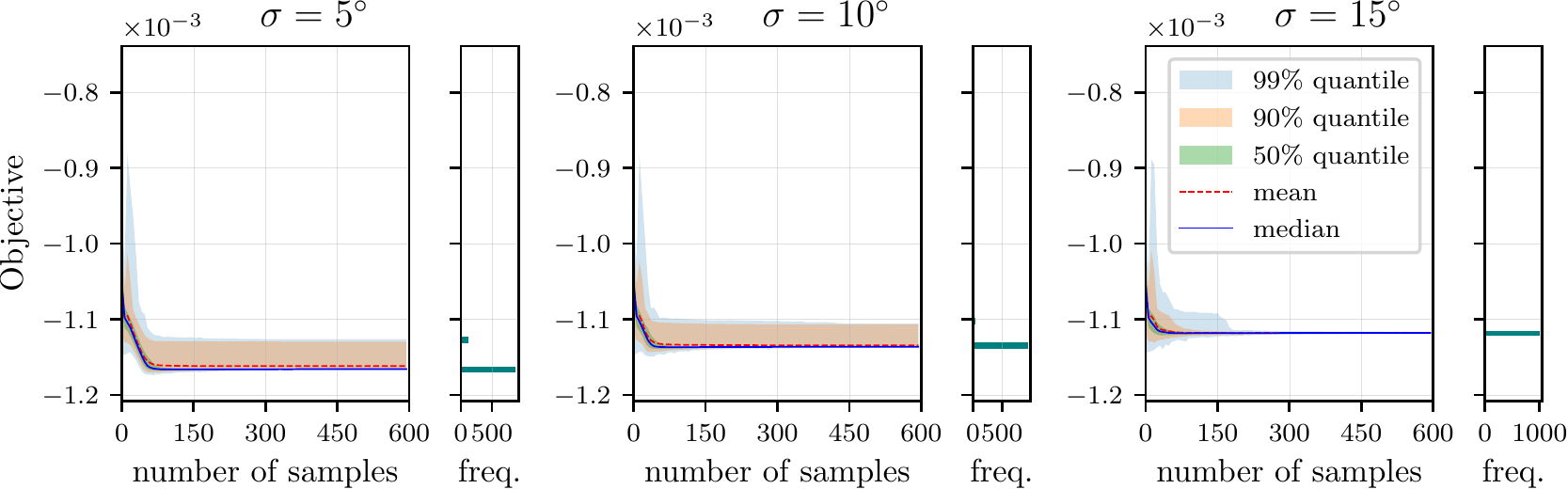}
	\end{center}
	\caption{
		Statistics of \num{1000} optimization runs with standard deviations $\sigma=\SI{5}{\degree}$, $\sigma=\SI{10}{\degree}$ and $\sigma=\SI{15}{\degree}$ of the orientations of the inclusions.
		The histograms represent the distribution of the final objective values.
		Each optimization was performed with \num{100} iterations, generating \num{6} samples per iteration.
		(For more information about the color references in this figure, the reader is referred to the digital version of this article.)
	}
	\label{fig:objective_vs_samples_pointcloud2d}
	\bigskip
	\includegraphics[width=\textwidth]{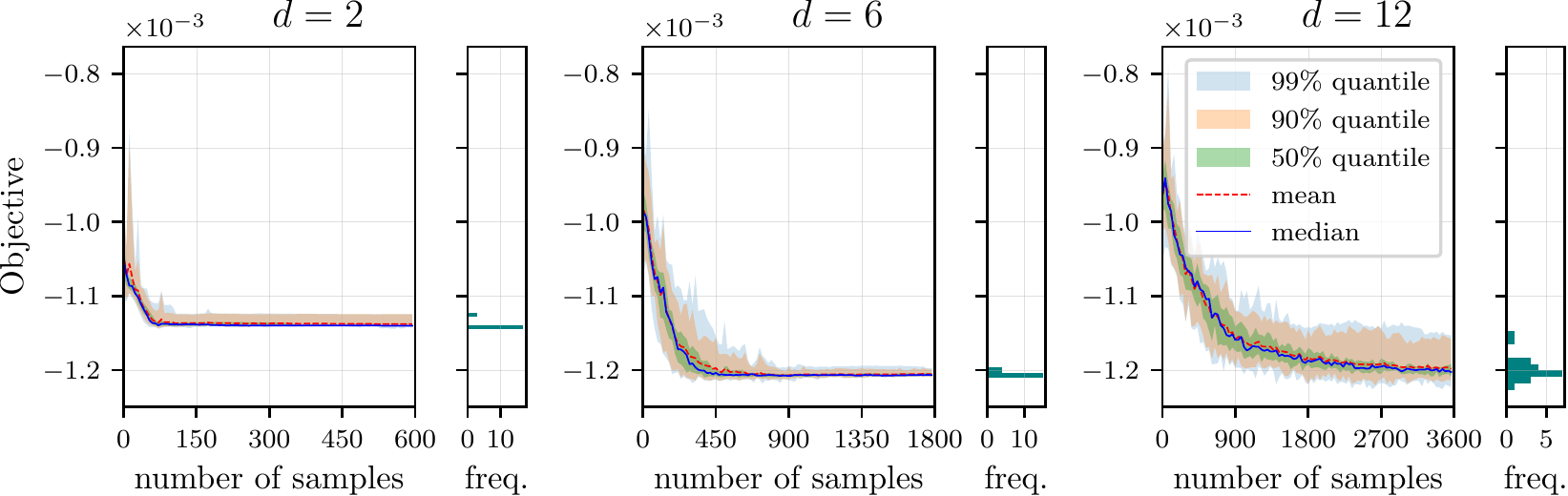}
	\caption{
		Statistics of \num{20} optimization runs with problem dimensions $d=\num{2}, \num{6}$ and $12$.
		The histograms represent the distribution of the final objective values.
		Each optimization was performed with \num{100} iterations, generating $3 \times d$ samples per iteration.
		The standard deviation parameters are $\sigma_\theta = \SI{3.6}{\degree}$ for orientations and $\sigma_r = \num{1.76e-3}$ for radii.
		(For more information about the color references in this figure, the reader is referred to the digital version of this article.)
	}
	\label{fig:objective_vs_samples}
\end{figure}

With a cheap-to-evaluate response surface approximation at hand, the optimization was run \num{1000} times, each with three different values of standard deviation: $\sigma=\SI{5}{\degree}, \sigma=\SI{10}{\degree}$ and $\sigma=\SI{15}{\degree}$, with $\sigma^\mathrm{max} = \SI{36}{\degree}$.
Due to the \SI{180}{\degree} periodicity of the orientations, the design variables are unbounded.
Fig.~\ref{fig:objective_landscape_2d_smoothed} shows a scatter plot of the final designs over the smooth objective function landscapes.
With increasing the degree of smoothing, the behavior of the objective function landscape becomes more global, implying less number of local minima.
Fig.~\ref{fig:objective_vs_samples_pointcloud2d} shows the statistics of the objective values with respect to the optimization iteration, which indicates that the final objective value, which is the average over a finite region in the design space, increases with standard deviation $\sigma$.

\subsection{Optimization with $2, 6$ and $12$ DoFs}

\begin{figure}[p]
	\centering
	\includegraphics[width=\textwidth]{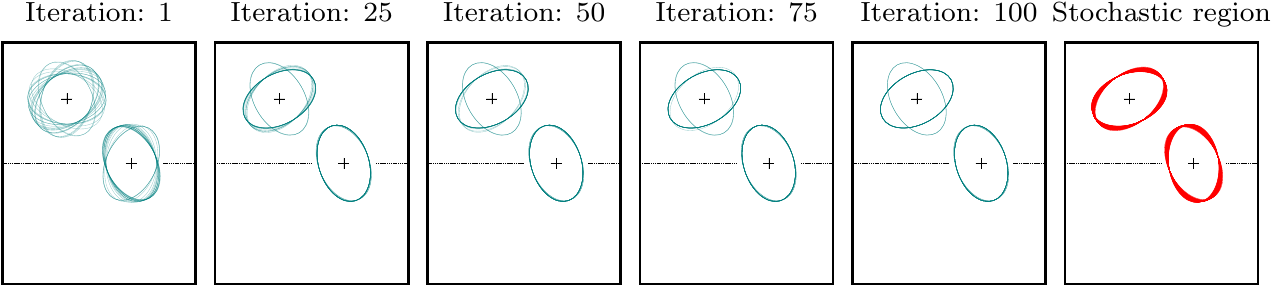}
	\caption{
		Design vs.\ optimization iteration for problem dimension $d=2$, gathered over \num{20} optimization runs.
		The stochastic region depicts random perturbations of the final design corresponding to the median of the final objective values from all the optimization runs.
		To better illustrate the independent design variables, only two inclusions are shown here.
		Please refer to Fig.~\ref{fig:design_parameters} for details on the design parameterization.
	}
	\label{fig:design_d2}
	\bigskip
	\includegraphics[width=\textwidth]{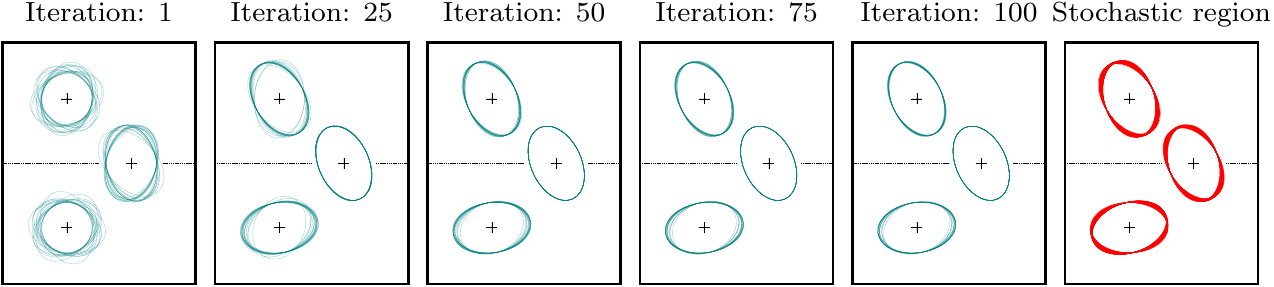}
	\caption{
		Design vs.\ optimization iteration for problem dimension $d=6$, gathered over \num{20} optimization runs.
		The stochastic region depicts random perturbations of the final design corresponding to the median of the final objective values from all the optimization runs.
		To better illustrate the independent design variables, only three inclusions are shown here.
		Please refer to Fig.~\ref{fig:design_parameters} for details on the design parameterization.
	}
	\label{fig:design_d6}
	\bigskip
	\includegraphics[width=0.82\textwidth]{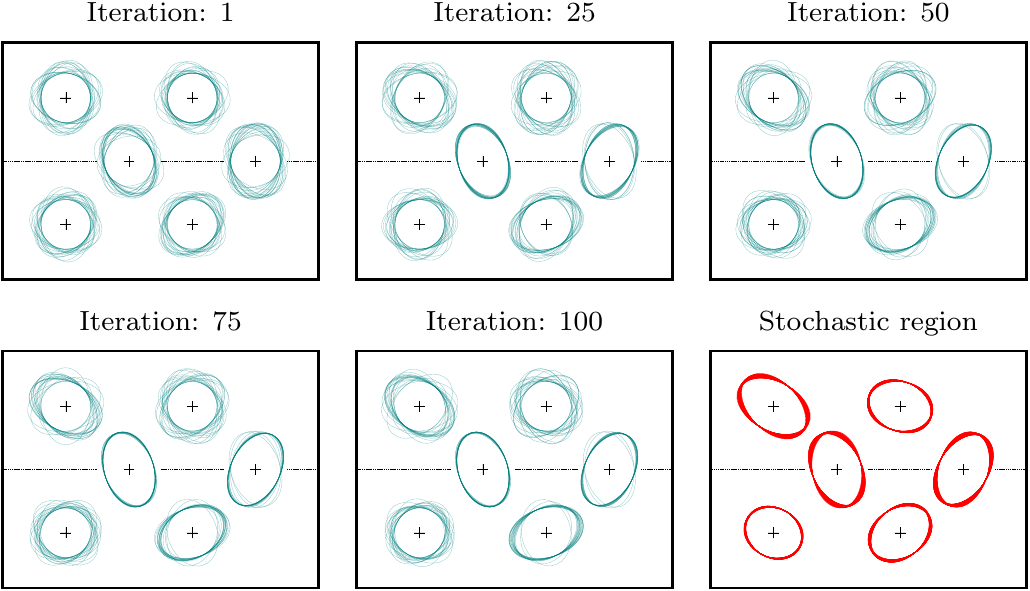}
	\caption{
		Design vs.\ optimization iteration for problem dimension $d=12$, gathered over \num{20} optimization runs.
		The stochastic region depicts random perturbations of the final design corresponding to the median of the final objective values from all the optimization runs.
		To better illustrate the independent design variables, only six inclusions are shown here.
		Please refer to Fig.~\ref{fig:design_parameters} for details on the design parameterization.
	}
	\label{fig:design_d12}
\end{figure}

Next, we test the optimization algorithm with three cases of problem dimension: $d=2, 6$ and $12$, for which the design parameters are described in Fig.~\ref{fig:design_parameters}.
Here, the design parameters are the orientations $\theta_i$ and the major radii $r_i$ of the elliptical inclusions, where the major radii are subjected to box constraints ($r_i \in [0.156, 0.244]$) to avoid possible overlap of the inclusions.
For the $d=2$ case, the major radii were fixed to the upper bound of \num{0.244}.
The minor radii of the inclusions are set to \num{0.15}.
The standard deviations were chosen to be $\sigma_\theta = \SI{3.6}{\degree}$ for the orientations and $\sigma_r = \num{1.76e-3}$ for the radii.
For the optimization runs, the design space is scaled such that all the design parameters use the same standard deviation parameter $\sigma$, and $\sigma^\mathrm{max}$ is chosen to be ten times larger.
The optimization procedure was performed \num{20} times for each of the three cases of the problem dimension.

\begin{figure}
	\begin{center}
		\includegraphics[width=\textwidth]{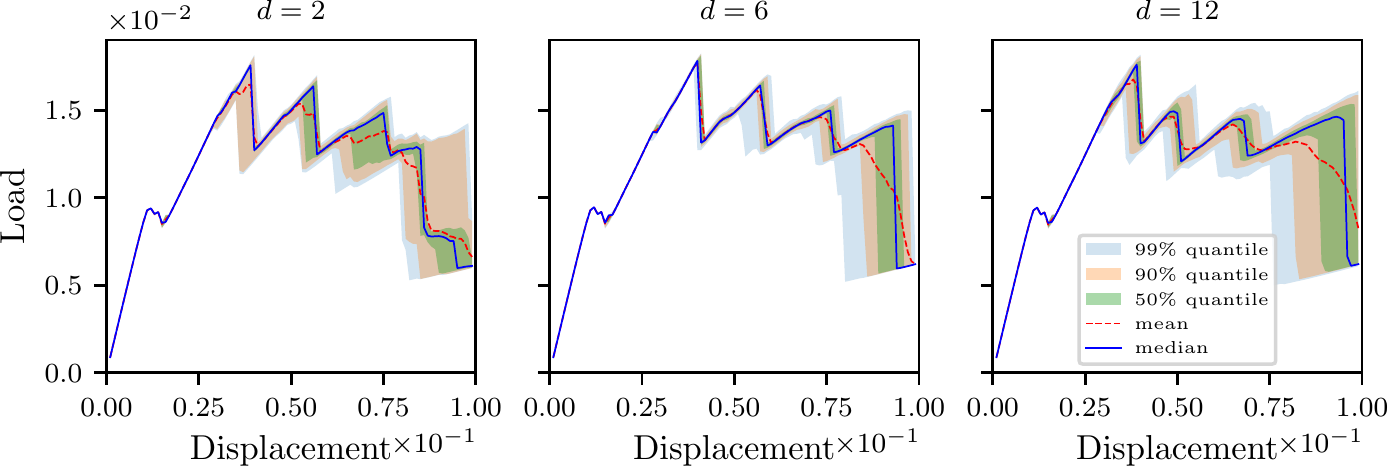}
	\end{center}
	\caption{
		Statistics of the load--displacement curves for the final design samples ($3 \times d$ samples per optimization run)  gathered over \num{20} optimization runs.
		(For more information about the color references in this figure, the reader is referred to the digital version of this article.)
	}
	\label{fig:load_disp_final}
\end{figure}

\begin{figure}
	\centering
	\begin{subfigure}{\textwidth}
		\centering
		\includegraphics[width=0.4\textwidth]{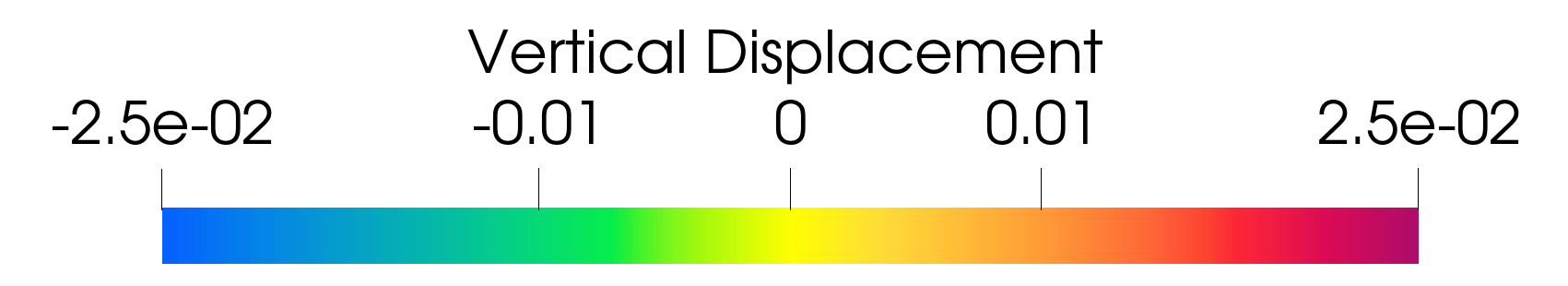}
	\end{subfigure}
	\begin{subfigure}{0.49\textwidth}
		\centering
		\fbox{\includegraphics[width=0.9\textwidth]{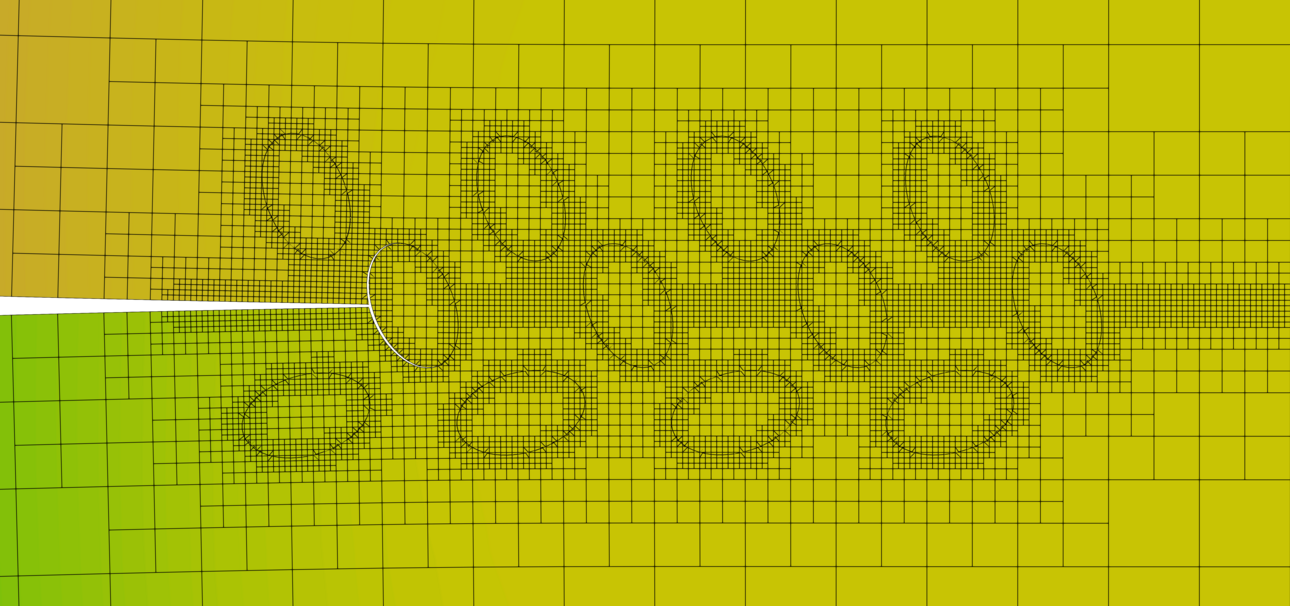}}
    \caption*{\normalsize{$t=0.020$}}
	\end{subfigure}
	\begin{subfigure}{0.49\textwidth}
		\centering
		\fbox{\includegraphics[width=0.9\textwidth]{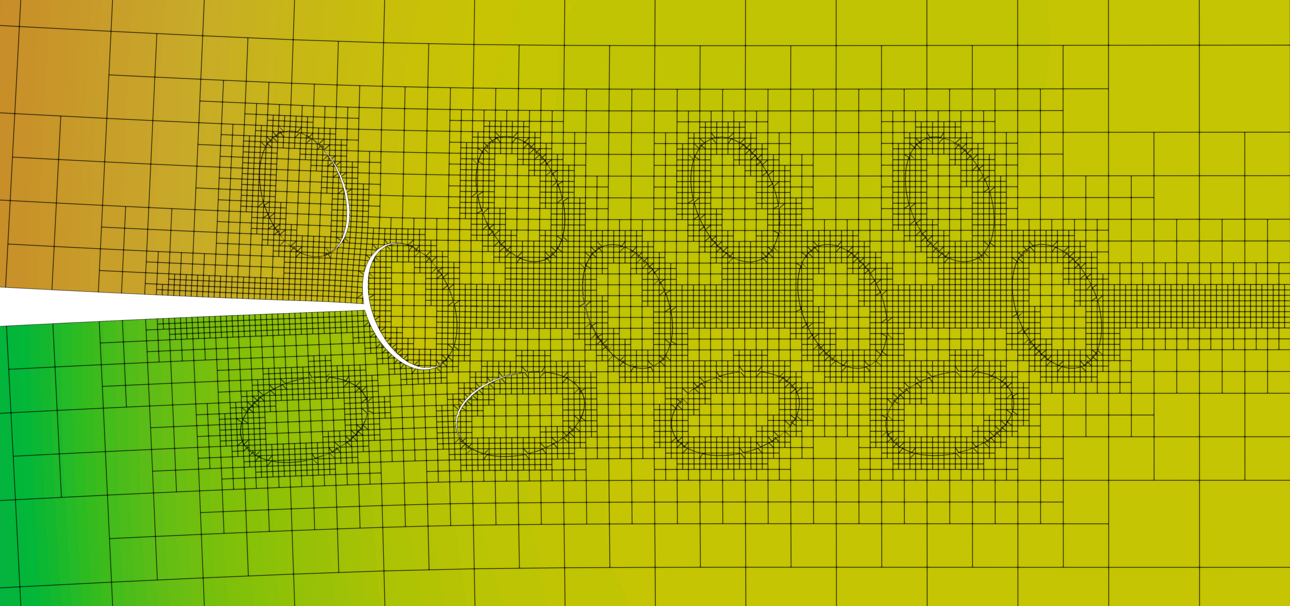}}
    \caption*{\normalsize{$t=0.040$}}
	\end{subfigure}
	\par\bigskip
	\begin{subfigure}{0.49\textwidth}
		\centering
		\fbox{\includegraphics[width=0.9\textwidth]{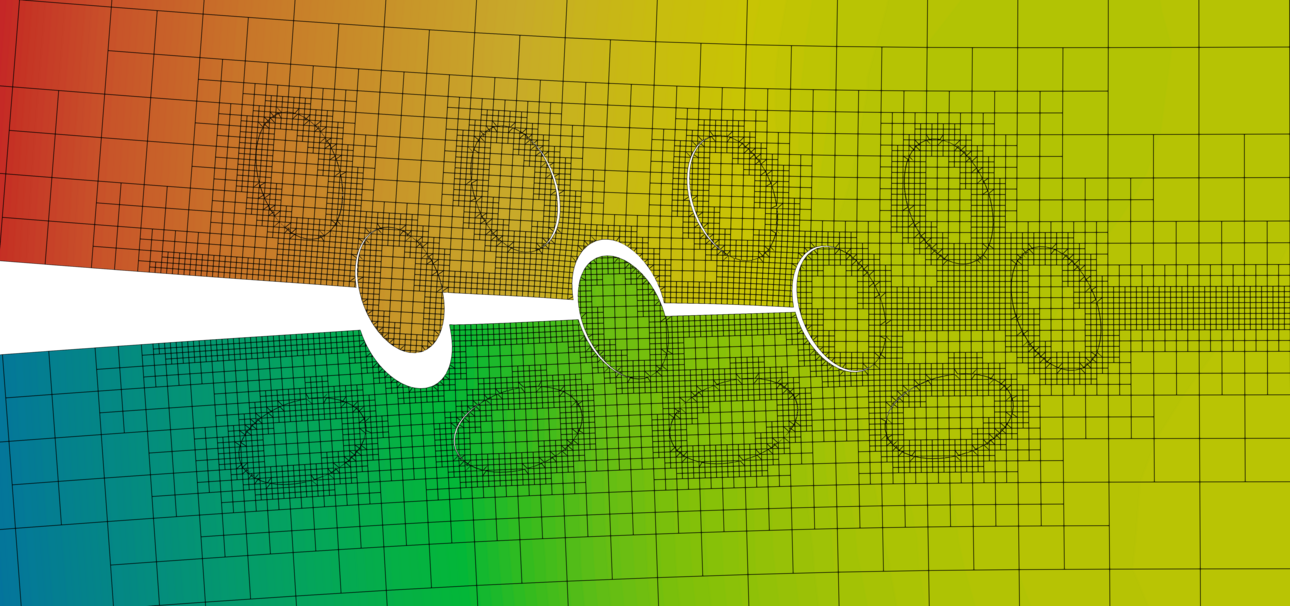}}
    \caption*{\normalsize{$t=0.060$}}
	\end{subfigure}
	\begin{subfigure}{0.49\textwidth}
		\centering
		\fbox{\includegraphics[width=0.9\textwidth]{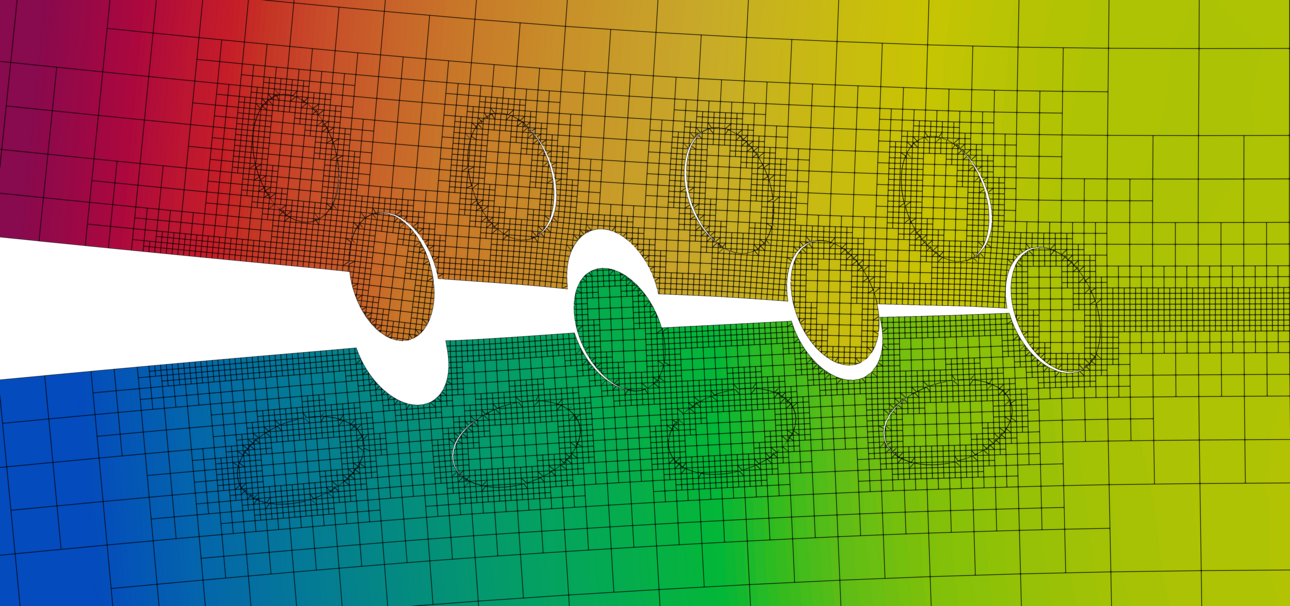}}
    \caption*{\normalsize{$t=0.080$}}
	\end{subfigure}
	\caption{
		Evolution of cracks for the final design obtained from optimization run with $d=6$ design parameters.
		The deformations are scaled by a factor of $10$ for better visualization of the crack openings.
		(For more information about the color references in this figure, the reader is referred to the digital version of this article.)
	}
	\label{fig:crack_evolution_visualization}
\end{figure}

\begin{figure}
	\centering
	\begin{subfigure}{0.44\textwidth}
		\centering
		\fbox{\includegraphics[width=0.9\textwidth]{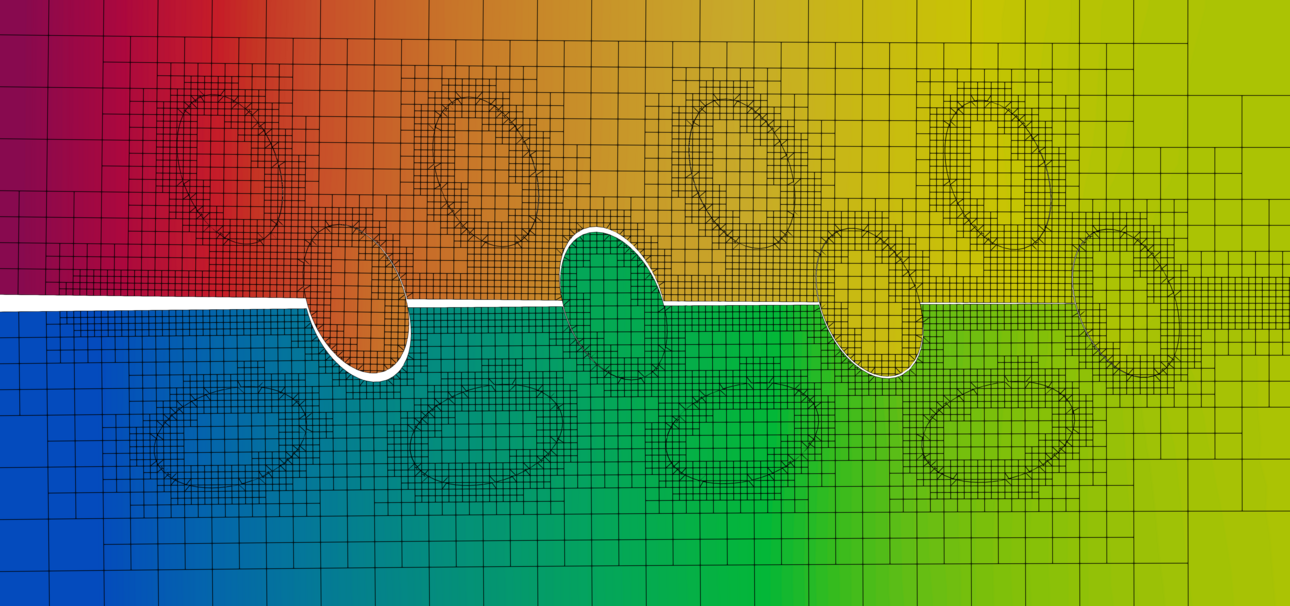}}
    \caption*{\normalsize{$t=0.094$}}
	\end{subfigure}
	\begin{subfigure}{0.44\textwidth}
		\centering
		\fbox{\includegraphics[width=0.9\textwidth]{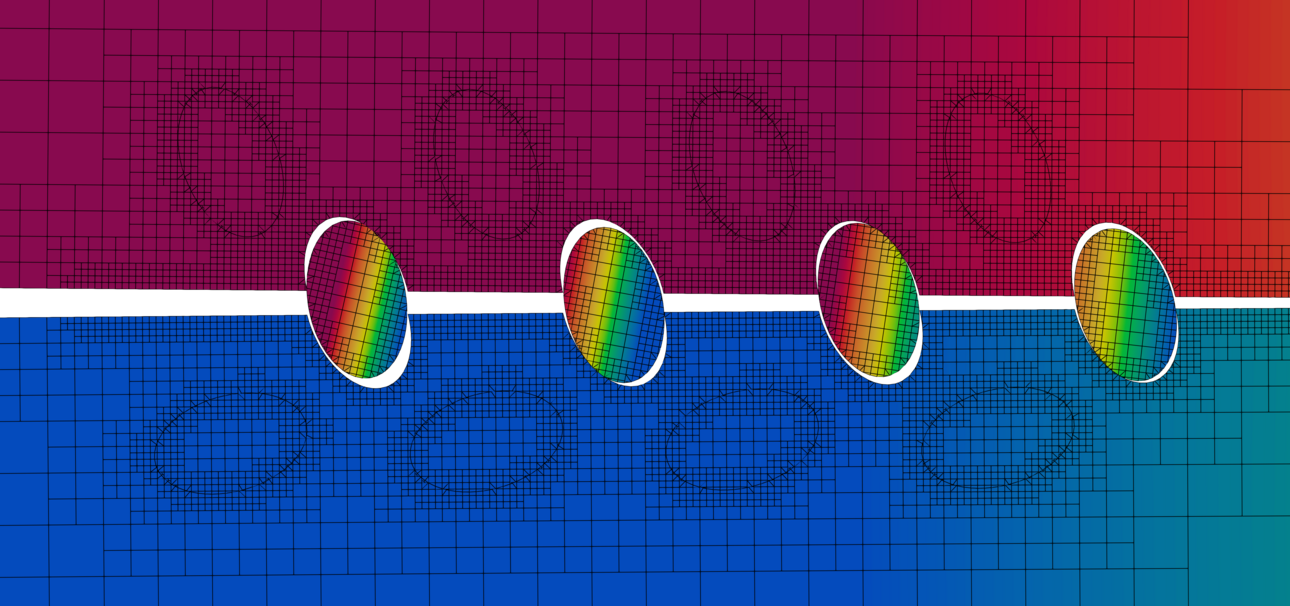}}
    \caption*{\normalsize{$t=0.095$}}
	\end{subfigure}
	\begin{subfigure}{0.1\textwidth}
		\centering
		\includegraphics[width=\textwidth]{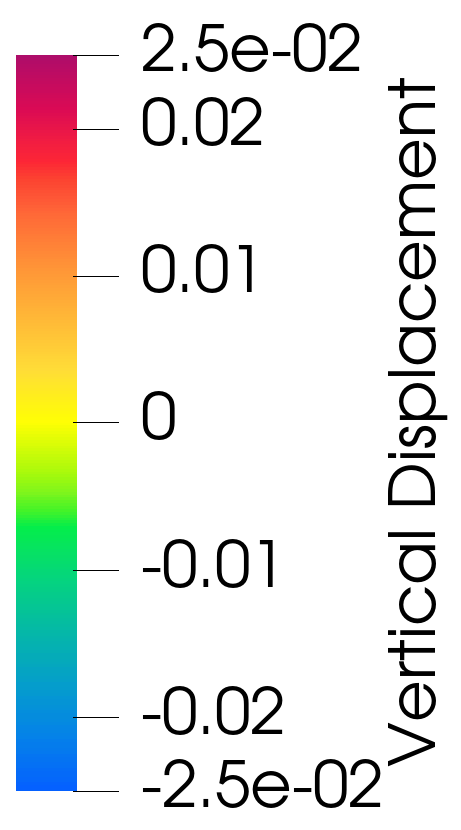}
		\caption*{}
	\end{subfigure}
	\caption{
		Fracture toughening as a result of crack bridging mechanism, illustrated for the final design obtained from optimization run with $d=6$ design parameters.
		(For more information about the color references in this figure, the reader is referred to the digital version of this article.)
	}
	\label{fig:crack_bridging}
\end{figure}

Fig.~\ref{fig:objective_vs_samples} shows the statistics of the objective value with respect to the optimization iteration.
There is a significant improvement ($\approx \SI{5.2}{\percent}$) in the objective value when the number of design DoFs is increased from \num{2} to \num{6}, while the shift to \num{12} DoFs does not lead to a further enhancement.
Moreover, there is a high variance in the objective value for $d=12$ case, which is explained by the fact that the nearest-neighbor interpolation is affected by the curse of dimensionality.
For high dimensions, a significantly larger number of samples is required to accurately predict the expected QoI and its gradients w.r.t.\ the design variables.
Fig.~\ref{fig:design_d2}, Fig.~\ref{fig:design_d6} and Fig.~\ref{fig:design_d12} show the final designs for all the optimization runs with $d=2, 6$ and $12$ design DoFs, respectively.
To account for the geometric symmetry of the arrangement of inclusions around the central horizontal axis of the beam, the designs are flipped to match the orientation of the leftmost central inclusion.
For $d=2$ and $d=6$, distinct local optima are observed, while for $d=12$, the low sensitivity of the objective w.r.t.\ some of the design parameters and the low accuracy of the nearest-neighbor interpolant in high dimensions lead to high variability of the final designs across all the optimization runs.
Fig.~\ref{fig:load_disp_final} shows load--displacement diagrams for the final designs.
Fig.~\ref{fig:crack_evolution_visualization} and Fig.~\ref{fig:crack_bridging} depict evolution of the cracks with the final design corresponding to the median objective value for the $d=6$ optimization case.
As shown by Fig.~\ref{fig:crack_bridging}, the improvement in the delamination resistance is due to the occurrence of crack bridging mechanism that provides additional resistance to crack growth.

\setlength{\fboxsep}{0pt}

\section{Conclusion}\label{sec:conclusion}

In this paper, a robust design optimization scheme was presented for enhancing delamination resistance of heterogeneous structures, for which the structural responses are discontinuous and exhibit numerical noise.
The method benefits from the response surface approximation via nearest-neighbor interpolation and its subsequent stochastic relaxation.
The use of smooth Gaussian distribution for the perturbations of the design variables allows computation of gradients to be used for searching local optima.
The proposed optimization algorithm was tested for design problems with $2$, $6$ and $12$ degrees of freedom.

In nutshell, the model offers the following advantages:
\begin{itemize}
	\item \emph{Robust design}.
	      In practical situations, it is crucial for the design to be robust under small deviations. This is achieved by incorporating a stochastic optimization approach where the QoI is optimized in an average sense.

	\item \emph{Efficiency}.
	      All data points from the previous optimization iterations are utilized for the evaluation of the expected QoI using the nearest-neighbor interpolation scheme.

	\item \emph{Derivative free}.
	      The optimization algorithm does not require derivatives of the underlying QoI with respect to the design variables.
	      This allows to incorporate fracture models of arbitrary complexity within the optimization framework.

	\item \emph{Parallel computation}.
	      The QoI evaluation for the proposed samples in each optimization iteration can be performed independently and in parallel on multiple compute resources.

	\item \emph{Avoiding premature convergence}.
	      Starting the optimization with a larger smoothing radius and gradually reducing it to the desired one helps overcome poor local optima and avoids premature convergence.
\end{itemize}

\textbf{Limitations}. Although, the optimization algorithm works well in low dimensions, the accuracy of the nearest-neighbor interpolation drops significantly in high dimensions, as observed for the \num{12} DoFs case.
Methods based on dimensionality reduction (e.g.\ using active subspaces~\cite{ActiveSubspaceConsta2015}) and distance-metric learning (e.g.~\cite{pmlr-v2-weinberger07a}) for the nearest-neighbor interpolation could improve the estimation of the QoI and thus the efficiency of the proposed stochastic optimization framework in high dimensions.

\section*{Acknowledgements}
The work was funded by the Deutsche Forschungsgemeinschaft (DFG, German Research Foundation) -- 377472739 / GRK 2423/1-2019.
The authors are grateful for the administrative and advisory support from the
Competence Unit for Scientific Computing (CSC)
at FAU Erlangen-N\"urnberg, Germany.


\printbibliography

\end{document}